\documentclass[9pt,a4paper,twoside]{rho-class/rho}
\usepackage[english]{babel}

\usepackage{amsmath,amssymb,amsthm}

\usepackage{blindtext}
\usepackage{algorithm2e}
\usepackage{cite}
\usepackage{latexsym}
\usepackage{indentfirst}
\usepackage{subfigure}
\usepackage[graphicx]{realboxes}
\usepackage{float}
\usepackage{graphicx}
\usepackage{placeins}
\usepackage{booktabs}
\usepackage{algorithmic}
\usepackage{multirow}
\usepackage{xcolor}
\allowdisplaybreaks[4]
\theoremstyle{plain}
\newtheorem{theorem}{Theorem}
\newtheorem{lemma}[theorem]{Lemma}

\theoremstyle{definition}

\newtheorem*{fact}{Fact}
\newtheorem*{definition}{Definition}
\SetKwComment{Comment}{/* }{ */}
\theoremstyle{remark}
\newtheorem{remark}{Remark}

\DeclareMathOperator{\argmin}{argmin}

\usepackage{tikz}

\newcommand{\bbm}{\begin{bmatrix}}
\newcommand{\ebm}{\end{bmatrix}}

\newcommand{\A}{\mathcal{A}}

\definecolor{MSBblue}{rgb}{0.204, 0.353, 0.541}
\definecolor{MSLightBlue}{rgb}{0.31, 0.506, 0.741}
\definecolor{Myred}{rgb}{0.71, 0.22, 0.29}

\newcommand{\yplu}[1]{{\color{orange}{ [\textbf{Yiping:}} #1]}}

\usepackage{tcolorbox}

\tcbuselibrary{breakable}

\usepackage{lastpage}
\allowdisplaybreaks[4]

\usepackage{lastpage}
\newtcolorbox{boxH}{
breakable,
    boxrule = 0pt, 
    leftrule = 6pt 
}

\newtcolorbox{boxK}{
breakable,
    sharpish corners, 
    boxrule = 0pt,
    toprule = 4.5pt, 
}

\usepackage{xr}
\makeatletter
\newcommand*{\addFileDependency}[1]{
\typeout{(#1)}
\@addtofilelist{#1}
\IfFileExists{#1}{}{\typeout{No file #1.}}
}\makeatother

\newcommand*{\myexternaldocument}[1]{%
\externaldocument{#1}%
\addFileDependency{#1.tex}%
\addFileDependency{#1.aux}%
}

\myexternaldocument{supplement}

\newtcbox{\alertinline}[1][red]
  {on line, arc = 0pt, outer arc = 0pt,
    colback = #1!20!white, colframe = #1!50!black,
    boxsep = 0pt, left = 1pt, right = 1pt, top = 2pt, bottom = 2pt,
    boxrule = 0pt, bottomrule = 1pt, toprule = 1pt}
\setbool{rho-abstract}{true} 
\setbool{corres-info}{true} 

\journalname{arXiv preprint. Under review.}
\title{Randomized Iterative Solver as Iterative Refinement:\\A Simple Fix Towards Backward Stability}

\author[1]{Ruihan Xu}
\author[2]{Yiping Lu}

\affil[1]{1Department of Mathematics, University of Chicago}
\affil[2]{2
Industrial Engineering \& Mangenement Science, Northwestern University}

\dates{This manuscript was compile on Oct 1, 2024}

\leadauthor{Xu et al.}
\footinfo{Creative Commons CC BY 4.0}
\smalltitle{Sketched Iterative and Recursive Refinement}
\institution{Northwestern University}
\theday{May 21, 2024} 

\corres{Yiping Lu}
\email{yiping.lu@northwestern.edu}

\license{\ \ This document is licensed under Creative Commons CC BY 4.0.\ccLogo}

\begin{abstract}
Iterative sketching and sketch-and-precondition are well-established randomized algorithms for solving large-scale over-determined linear least-squares problems. In this paper, we introduce a new perspective that interpreting Iterative Sketching and Sketching-and-Precondition as forms of Iterative Refinement. We also examine the numerical stability of two distinct refinement strategies: iterative refinement and recursive refinement, which progressively improve the accuracy of a sketched linear solver. Building on this insight, we propose a novel algorithm, Sketched Iterative and Recursive Refinement (SIRR), which combines both refinement methods. SIRR demonstrates a \emph{four order of magnitude improvement} in backward error compared to iterative sketching, achieved simply by reorganizing the computational order, ensuring that the computed solution exactly solves a modified least-squares system where the coefficient matrix deviates only slightly from the original matrix. To the best of our knowledge, \emph{SIRR is the first asymptotically fast, single-stage randomized least-squares solver that achieves both forward and backward stability}.
\end{abstract}


\begin{document}
	
    \maketitle
    \thispagestyle{firststyle}


\section{Introduction}

\rhostart{R}andomized numerical linear algebra (RNLA) \cite{tygert2009fast,rokhlin2008fast,woodruff2014sketching, martinsson2020randomized, drineas2016randnla, mahoney2011randomized, avron2010blendenpik} is a rapidly evolving branch of matrix computations, driving significant progress in low-rank approximations, iterative methods, and projections. This field has demonstrated that randomized algorithms are highly effective tools for developing approximate matrix factorizations. These methods are remarkable for their simplicity and efficiency, often producing surprisingly accurate results.

In this paper, we consider randomized algorithms to solve the overdetermined linear least-squares problem 
{
\begin{equation}
   x = \arg\min_{y \in \mathbb{R}^n} \|b - Ay\| \quad (A \in \mathbb{R}^{m \times n}, b \in \mathbb{R}^m)
\end{equation}
}
where \(\|\cdot\|\) denotes the Euclidean norm. This is one of the core problems in computational sceience \cite{lim2017fast,fitzsimons2016improved,lin2016approximating,tropp2019streaming}, {statistics \cite{yang2017randomized,khoo2024nonparametric} and accelerating machine learning \cite{avron2017faster,ahle2020oblivious,hur2023generative}}.  In the past two decades, researchers in the field of \textit{randomized numerical linear algebra} \cite{martinsson2020randomized,mahoney2011randomized,woodruff2014sketching} have developed least-squares solvers that are faster than Householder QR factorization \cite{golub2013matrix}, the textbook algorithm for least square, which runs in \(O(m n^2)\) operations. Randomized algorithms first sketch $A$ to a smaller matrix $SA$ with a random sketch matrix $S\in\mathbb{R}^{[cn]\times m}$ for some constant $c>1$. The random embedding $v\rightarrow Sv$ satisfies $\|Sv\|\approx \|v\|$ for all vectors $v\in\text{range}([A b])$ and matrix–vector products $v\rightarrow Sv$ can be computed efficiently \cite{mahoney2011randomized,martinsson2020randomized}. 

There are two main approaches to using the sketched matrix $SA$ for a fast randomized least squares solver: \emph{the sketch-and-precondition} \cite{rokhlin2008fast} method and \emph{iterative Hessian sketching} \cite{pilanci2016iterative,wang2017sketching,ozaslan2019iterative}. Most of the solvers (e.g. Blendenpik \cite{avron2010blendenpik}) have a complexity of $O(mn\log m)$ operations. This is significantly better than the $O(mn^2)$ complexity. Consequently, for large least-squares (LS) problems, randomized solvers can be substantially faster than the LS solver implemented in LAPACK \cite{avron2010blendenpik}. However, recent research \cite{meier2024sketch,epperly2023fast} surprisingly finds that sketch-and-precondition \cite{rokhlin2008fast,avron2010blendenpik} and iteratively Hessian sketch \cite{pilanci2016iterative,wang2017sketching,ozaslan2019iterative} are numerically unstable in their standard form, both stagnate in terms of residual and backward error, potentially before optimal levels are reached. \cite{meier2024sketch} further propose sketch-and-apply, which is a provable method that attains backward stable solutions under modest conditions.  Unfortunately, sketch-and-apply requires $O(mn^2)$ operations, the same as Householder QR-based direct solvers. In this paper, we provide a definitive answer to the open question posed by \cite{meier2024sketch, epperly2023fast}:

{
\begin{center}
    {\emph{Is there a randomized least-squares algorithm that is both (asymptotically) faster
than Householder QR and numerically stable?}}
\end{center}}

We constructed a solver called \textbf{S}ketched  \textbf{I}terative and \textbf{R}ecursive \textbf{D}ebiasing, which enjoys both forward and backward stability while requires only $O(mn+n^3)$ computation. Our approach is based on a novel, unified perspective on sketch-and-precondition methods and iterative Hessian sketching. Although these two techniques may seem different, we demonstrate that they can be interpreted as iterative refinement processes. Iterative refinement (IR) is a well-known method for solving linear systems by progressively improving the accuracy of an initial approximation. We show that employing iterative refinement, a sketch-and-solve solver is equivalent to using Jacobi iteration in a sketch-and-precondition framework. We investigated the conditions that a single-step approximate solver needs to satisfy in order for iterative refinement to potentially achieve backward stability. To construct the single-step approximate solver, we studied another way for iterative refinement called \textbf{S}ketched \textbf{R}ecursive \textbf{R}efinement.  Note that we find, both theoretically and numerically, that only in certain cases where data noise is relatively large, SRR alone can achieve a backward stable solution. Only using SRR as the meta-algorithm of iterative refinement, \emph{i.e.} \textbf{S}ketched  \textbf{I}terative and \textbf{R}ecursive \textbf{D}ebiasing, can provide a backward stable algorithm.

We would like to highlight a concurrent work \cite{epperly2024fast}, which also developed a backward stable solver with a computational complexity of $O(mn+n^3)$. However, the FOSSILS solver proposed in their work follows a two-stage approach, where each stage involves an iterative process. In contrast, our algorithm is a single-stage solver that offers the flexibility to stop at any point during the computation, making it more adaptable for scenarios where early termination is necessary or beneficial.

\paragraph{Notation} Through out this paper, $A\in \mathbb{R}^{m\times n}$, $S\in\mathbb{R}^{s\times m}$, $b\in \mathbb{R}^m$. $\|\cdot\|$ denotes vector $\ell_2$ norm for vectors and operator $\ell_2$ norm for matrices. We use $|\cdot|$ denotes $l_1$ norm. $\kappa = \|A\|\|A^\dag\|$ is the condition number of $A$ and $\sigma_{\max}(\cdot)$,$\sigma_{\min}(\cdot)$ denotes the largest and smallest singular value. $u$ denotes the machine epsilon which is used to measure the level of roundoff error in the floating-point number system. For IEEE standard double precision, $u$ is around $2\times 10^{-16}$. $a \lesssim b$ denotes $a\leq cb$ for some small constant $c$, which is independent of $m,n,s,\kappa,u$. $a \asymp b$ indicates that $a\lesssim b$ and $b\lesssim a$. $\gamma_m$ is defined as $\gamma_m=\frac{mu}{1-mu}$. In numerical analysis, we assume that $u\kappa n^{\frac{3}{2}}<1$ and $un^\frac{3}{2}\|x\ast\|\lesssim \|b-Ax^\ast\|$. We also assume that $m$ has the same order with $n$ for computational simplicity, which will be restated in the following sections. Note that $n^{\frac{3}{2}}<1$ is a guarantee for a nonsingular $\hat{R}$ computed in QR factorization according to \cite[Theorem 19.3]{higham2002accuracy}. Without loss of generality, $\|A\|=\|b\|=1$ is assumed in analysis, except in forward stability analysis where we keep $\|A\|$ and $\|b\|$ unknown to align with Wedin's perturbation theorem. Computed quantities wear a hat, e.g. $\hat{x}$ denotes the computed approximation of $x$.

\subsection{Contribution}

We offer a unified understanding of existing randomized least squares solvers, such as iterative sketching and sketch-and-precondition, by interpreting them as forms of iterative refinement. This new perspective enables the development of novel techniques for analyzing the numerical stability of randomized algorithms by exploring and comparing the stability of iterative and recursive refinement strategies for progressively improving the accuracy of sketched linear solvers. Based on the analysis, we propose Sketched Iterative and Recursive Refinement (SIRR), which combines iterative and recursive refinement techniques and achieves the first single stage provably backward stable and computationally efficient, with asymptotic complexity $O(mn+n^3)$, faster than traditional direct solvers.

\section{Preliminary}

\paragraph{Sketch-and-Precondition} There are lots of randomized methods that obtain a right preconditioner from $SA$ for further iterative LS method,  which is known as sketch-and-precondition \cite{rokhlin2008fast, avron2010blendenpik, meng2014lsrn}. The core insight of Sketch-and-Precondition is that sketching matrices can be used to precondition (i.e., reduce the condition number) the original matrix $A\in\mathbb{R}^{m\times n}$. To be specific, for a matrix $A\in\mathbb{R}^{m\times n}$ and sketching matrix $S\in\mathbb{R}^{s\times m}$ with distortion $0<\eta<1$ (\emph{i.e.} $(1-\eta)\|Ay\|\le \|SAy\|\le (1+\eta)\|Ay\|$ holds for all $y\in\mathbb{R}^n$), the preconditioner $R$ can be obtained from QR factorization of matrix $SA=QR$ with $Q$ orthonormal and $R$ square. The preconditioner $R$ satisfies
$$
\frac{1}{1+\eta}\le \sigma_{\min}(AR^{-1})\le \sigma_{\max}(AR^{-1})\le \frac{1}{1-\eta}.
$$
To be specific, one can always construct a random sparse embedding matrix $S$ that satisfies the following Lemma.

\begin{lemma}[\cite{higham2002accuracy,meier2024sketch,epperly2023fast, epperly2024fast}]
\label{numerical error}
        For matrix $A\in \mathbb{R}^{m\times n}$, there exists sketching matrix $S\in \mathbb{R}^{s\times m}$. Suppose that $\hat{R}\hat Q=SA$ is the QR decomposition of matrix $SA$, then the following inequalities holds:
        \begin{itemize}
\setlength{\itemsep}{0pt}
\setlength{\parsep}{0pt}
\setlength{\parskip}{0pt}
            \item $\|\hat{R}\|\lesssim \|A\|$, $\|\hat{R}^{-1}\|\lesssim \frac{\kappa}{\|A\|}$
            \item $1-u\kappa n^\frac{5}{2}\lesssim \sigma_{min}(A\hat{R}^{-1})\lesssim \sigma_{max}(A\hat{R}^{-1})\lesssim 1+u\kappa n^\frac{5}{2}$
        \end{itemize}
    \end{lemma}

One of the most prominent sketch-and-precondition techniques is using $R$ as the preconditoner for LSQR \cite{paige1982lsqr}  which is known as Blendenpik \cite{avron2010blendenpik}. In exact arithmetic, Blendenpik has a complexity of $\mathcal{O}(mn \log m)$ operations, which is better than the $\mathcal{O}(mn^2)$ QR-based direct solver. Consequently, for large LS problems, Blendenpik can be substantially faster than the LS solver implemented in LAPACK, a widely used software library for numerical linear algebra.

\paragraph{Iterative Hessian Sketching \cite{pilanci2016iterative,ozaslan2019iterative}} Iterative Sketching start from an initial solution $x_0\in\mathbb{R}^n$ generate iterates $x_1,x_2,\cdots$ by solving a sequence of the sketched least-squares problems
\begin{equation}
\label{eq:ihs} x_{i+1}=x_i+\argmin_{y\in\mathbb{R}^n}\frac{1}{2}\|(SA)y\|^2-y^\top A^\top (b-Ax_i),
\end{equation}
for $i=0,1,2,\cdots$. As with the classical least-squares sketch, the quadratic form is defined by the matrix $SA\in\mathbb{R}^{m\times d}$, which leads to computational savings. The closed form solution of (\ref{eq:ihs}) is given via $x_{i+1}=x_i+(A^\top S^\top SA)^{\dagger}A^\top (b-Ax_i)$ which encounter with the iterative refine a sketch-and-apply solver which shown in Algorithm \ref{alg:Iterativerefinement}. In Section \ref{section:equal}, we also show that Iterative Hessian Sketching/Iterative Refinement is equivalent to Sketch-and-Precondition using a Jacobi Iteration Solver.

\paragraph{Backward-Stability} 

Backward stability refers to the property of a numerical algorithm where the computed solution is the exact solution to a slightly perturbed version of the original problem. Specifically, a solver is said to be backward stable if the solver satisfies the following property:
{
\begin{definition}[Backward error]
    In floating point arithmetic, it produces a numerical solution $\hat{x}$ that is the exact solution to a slightly modified problem:

\begin{equation}
\hat{x} = \arg\min_{y \in \mathbb{R}^n} \|(b + \Delta b) - (A + \Delta A)y\|
\end{equation}

where the (relative) size of the perturbations is at most

\begin{equation}
\|\Delta A\| \leq c\|A\|, \quad \|\Delta b\| \leq c\|b\| \quad \text{provided } c < 1.
\end{equation} 
\end{definition}
}

\cite{epperly2024fast}
show that a backward stable solver can achieve accurate estimation of each component of the solution and can enforce residual orthogonality, \emph{i.e.} the KKT condition of the least square problem that $A^\top(Ax-b)=0$. The classic Householder QR least-squares method is backward stable \cite[Ch. 20]{higham2002accuracy}. However, recent works \cite{meier2024sketch,epperly2023fast} showed that randomized sketching solver is not backward stable.

 To prove a solver is backward stable, we follow \cite{gratton2012accuracy,epperly2024fast} which utilize the Karlson-Wald\'en estimate $\widehat{\mathrm{BE}}_\theta(\widehat{\boldsymbol{x}}):=\frac{\theta}{\sqrt{1+\theta^2\|\widehat{\boldsymbol{x}}\|^2}}\left\|\left(\boldsymbol{A}^{\top} \boldsymbol{A}+\frac{\theta^2\|\boldsymbol{b}-\boldsymbol{A} \widehat{\boldsymbol{x}}\|^2}{1+\theta^2\|\widehat{\boldsymbol{x}}\|^2} \mathbf{I}\right)^{-1 / 2} \boldsymbol{A}^{\top}(\boldsymbol{b}-\boldsymbol{A} \widehat{\boldsymbol{x}})\right\|$ which can estimate the backward error up to a constant, \emph{i.e.} $\widehat{\mathrm{BE}}_\theta(\widehat{\boldsymbol{x}})\leq\mathrm{BE}_\theta(\widehat{\boldsymbol{x}})\leq\sqrt{2}\widehat{\mathrm{BE}}_\theta(\widehat{\boldsymbol{x}})$ \cite{gratton2012accuracy}. Given singular value decomposition $A=\sum_{i=1}^n \sigma_iu_iv_i^\top$, the Karlson-Wald\'en estimation indicates that a least square $\hat x$ is backward stable is equivalent to satisfying a component-wise error bound $\left| \mathbf{v}_i^\top (\hat{\mathbf{x}} - \mathbf{x}) \right| \lesssim \sigma_i^{-1} \cdot (1 + \|\hat{\mathbf{x}}\|) u + \sigma_i^{-2} \cdot \| \mathbf{b} - A \hat{\mathbf{x}} \| u \quad \text{for} \ i = 1, \dots, n.$.
{

In our paper, we adopt the proof technique from \cite{epperly2024fast}, utilizing the following decomposition, which we incorporate as $\alpha-\beta$ accuracy in our manuscript.
    \begin{definition}[$\alpha-\beta$ Accuracy]{\small We define $\hat x$ is $\alpha-\beta$ accurate if there exists $e_1,e_2\in\mathbb{R}^n$ such that $\|e_1\|,\|e_2\|\leq 1$ and
    \begin{align*}
        \hat{x}-x^\ast=\alpha(1+\|\hat{x}\|)\hat{R}^{-1}e_1+\beta \|b-A\hat{x}\|(A^\top A)^{-1}e_2,
    \end{align*}
    {where $\hat{R}$ is a preconditioner of $A$ such that for any singular value of $A\hat R^{-1}$ satisfies $\sigma(A\hat R^{-1})\asymp 1$.}}
    \end{definition}
}
\begin{lemma}
\label{becondition}
    The computed solution $\hat{x}$ of problem $Ax=b$ has backward error $be(\hat{x})\lesssim\sqrt{n}\epsilon$ if 
    \begin{align}
    \label{beform}
        \hat{x}-x^\ast=\epsilon(1+\|\hat{x}\|)\hat{R}^{-1}e_1+\epsilon \|b-A\hat{x}\|(A^\top A)^{-1}e_2,
    \end{align}
    where $e_i\in\mathbb{R}^n$  satisfies $\|e_i\|\lesssim 1 (i=1,2)$.  
\end{lemma}

\paragraph{Numerical Stability} We provide several basic facts about numerical errors generated in floating-point arithmetic, most of which can be found in~\cite{higham2002accuracy}. For error analysis, we denote the numerical error of an expression computed in floating-point arithmetic as $\text{err}(\cdot)$. Specifically, for a real number $x$, let $\text{fl}(x)$ denote its floating-point approximation. The numerical error in $x$ is then defined as $\text{err}(x) = |x - \text{fl}(x)|$. {Recall that $u$ denotes the unit roundoff}, which is the maximum relative error in representing a real number in floating-point arithmetic. That is, for any real number $x$, we have $|\text{fl}(x) - x| \leq u |x|$. We also define $\gamma_n$ for a positive integer $n$ as
$\gamma_n = \frac{n u}{1 - n u}$,
assuming $n u \ll 1$, so that $\gamma_n \approx n u$.

\begin{fact}
    For vector $x,y\in \mathbb{R}^n$, matrix $A\in \mathbb{R}^{m\times n}$, upper triangular matrix $R\in \mathbb{R}^{n\times n}$, we have
    \begin{itemize}
\setlength{\itemsep}{0pt}
\setlength{\parsep}{0pt}
\setlength{\parskip}{0pt}
        \item $\|\text{err}(x\pm y)\|\leq u\sqrt{n}\|x\pm y\|$.
        \item $\|\text{err}(Ax)\|\leq \sqrt{n}\gamma_n\|A\|\|x\|$.
        \item {higham2002accuracy}For problem $Rx=y$, the solution by Gaussian-elimination satisfies $(R+E)(R^{-1}y+error(R^{-1}y))=y$ where $|E|\lesssim \gamma_n|R|$. This result further leads to $err(R^{-1}y) = \sqrt{n}\gamma_n\|R^{-1}y\|R^{-1}e$, where $\|e\|\lesssim 1$\cite[Theorem 8.5].
        \item For problem $Ax=b$, the solution by QR factorization satisfies $(A+\delta A)(A^\dag b + \text{err}(A^\dag b))=b+\delta b$, where $\|\delta A\|\lesssim \gamma_{n^2}\|A\|, \|\delta b\|\lesssim \gamma_{n^2}\|b\|$ \cite[Theorem 19.5]{higham2002accuracy}. 
    \end{itemize}
\end{fact}

\section{Randomized Solver As Iterative Refinement}

In this section, we present a novel approach for constructing a fast and stable randomized least squares solver by iteratively refining an approximate solver which we call a meta-algorithm, \emph{e.g.} sketch-and-apply or early stopped iterative randomized solver. We introduce two ways to do the refinement: iterative refinement and recursive refinement. Both refinement process starts from a meta-algorithm and improve the previous solution by correcting it based on the residual error. The key difference between iterative and recursive refinement processes is that iterative refinement improves the solution by applying the meta-algorithm at each step to correct the residual, while recursive refinement refines the solution by repeatedly applying the same current solver to the residual error.

\begin{algorithm}

\caption{\textbf{S}ketched \textbf{I}terative \textbf{R}efinement}

\label{alg:Iterativerefinement}
    \begin{tcolorbox}
    \hrule\hrule
    \textbf{SIR}: Sketched Iterative Refinement
    \hrule
    \SetKwInOut{Input}{Input}
    \SetKwInOut{Output}{Output}
    \SetKw{KwBy}{by}
    \SetKw{KwReturn}{Return}
    \SetKw{KwVia}{Via}
    \Input{1}
    \Output{2}
    \hrule\hrule
    
    \If{$N=0$}{\KwReturn $\text{SIR}^{\text{meta}}_0(b)$ \KwVia a meta-algorithm $\text{SIR}^{\text{meta}}_0(b)=\text{ALG}^{\text{meta}} (A^\top b)$\;
    
    \Comment*[r]{Initialization via Meta-Algorithm}}
    \For{$i\gets1$ \KwTo $N$ \KwBy $1$}{
    $\text{SIR}^{\text{meta}}_i(b):=\text{SIR}^{\text{meta}}_{i-1}(b)+{\color{orange}\text{ALG}^{\text{meta}}} (A^\top (b-A \cdot \text{SIR}^{\text{meta}}_{i-1}(b)))$\Comment*[r]{Iterative Refinement via Meta-Algorithm}
    }
    \KwReturn $\text{SIR}^{\text{meta}}_N(b)$
    \end{tcolorbox}
    \end{algorithm}
\subsection{Iterative and Recursive refinement}
\label{section:IRandRR}

\paragraph{Iterative Refinement} Iterative refinement \cite{wilkinson2023rounding,moler1967iterative,wu2023stable} is the classical approach to improving the quality of a computed solution in numerical linear
algebra. The idea of iterative refinement is simple, to improve the quality of an approximate solution $x_i$, solve for the error $\delta x_i=x-x_i$ via approximately solving $\delta x_i:=\arg\min_{\delta x_i}\|b-Ax_i-A\delta x_i\|$. Classically, the inexact solve used in the refinement step is a classical direct solver such as QR factorization computed in lower numerical precision (i.e., single precision), and all the other steps are performed in higher precision (e.g., double precision) \cite{golub1966note,carson2024comparison}. In our paper, we design an iterative algorithm, where each step incorporates the concept of iterative refinement, using a fast randomized linear solver to approximately solve the system. The algorithm is detailed in Algorithm \ref{alg:Iterativerefinement}.

\paragraph{Recursive Refinement} We also introduce a novel way to implement an iterative refinement process which we call it (sketched) recursive refinement approach. Sketched Recursive Refinement process also iteratively refines the solution by incorporating corrections from previous iterations. Different from iterative refinement which updates the current solution by applying a fixed procedure to adjust the solution, recursive refinement refers back to itself to perform the next step and solve the problem in a nested fashion. The algorithm is detailed in Algorithm \ref{alg:recursiverefinement}. Later, we demonstrate that recursive refinement is simply a reorganization of the computational steps in iterative refinement but the two types of refinement enjoy very different numerical stability behavior.

\begin{algorithm}
\caption{\textbf{S}ketched \textbf{R}ecursive \textbf{R}efinement.}\label{alg:recursiverefinement}

\begin{tcolorbox}
\hrule\hrule
\textbf{SRR}: Sketched Recursive Refinement
\hrule
\SetKwInOut{Input}{Input}
\SetKwInOut{Output}{Output}
\SetKw{KwBy}{by}
\SetKw{KwReturn}{Return}
\SetKw{KwVia}{Via}
\Input{1}
\Output{2}
\hrule\hrule

\If{$N=0$}{\KwReturn $\text{SRR}_0(b)$ \KwVia meta-algorithm $\text{ALG}^{meta}(A^\top b)$\;}
\For{$i\gets1$ \KwTo $N$ \KwBy $1$}{
    $\text{SRR}_i(b):=\text{SRR}_{i-1}(b)+{\color{orange}\text{SRR}_{i-1}}(A^\top b-A^\top  A\cdot \text{SRR}_{i-1}(b))$\Comment*[r]{Recursive Refinement}
    }
\KwReturn $\text{SRR}_N(b)$
\end{tcolorbox}
\end{algorithm}

\paragraph{Recursive Refinement as Reorganizing Computation} We would like to point out that Recursive refinement and Iterative refinement perform the same if one uses exact arithmetic. With a linear meta-algorithm, \emph{i.e.} $ALG^{meta}(A^\top b)$ can be represented as $TA^\top b+q$ for some matrix $T$ which includes most useful randomized solver such as Sketch-and-Apply, the results of $\text{SIR}_{N}(b)$ and $\text{SRR}_{\log_2 N}(b)$ are the same and both can be presented in the same form as geometric series as $x=\sum_{i=0}^{N} (I-TA)^{i}Tb$ with same amount of compute $O(Nmn)$. This means that Recursive Refinement is just a reorganization of computation order in the Iterative Refinement procedure and would generate the same computational result if one use exact arithmetic.  However, in the following discussion, we show that Recursive Refinement and Iterative Refinement behave very differently when using a floating point arithmetic.

\paragraph{Equivalence between Iterative Refinement and Sketch-and-Precondition}
\label{section:equal}

Iterative Refinement (Iterative Hessian Sketching) and the Sketch-and-Precondition approach are commonly regarded as two distinct methodologies for designing iterative randomized least squares solvers. In this remark, we demonstrate the surprising equivalence between sketched iterative refinement and the sketch-and-precondition method. This insight provides a unified perspective on modern randomized linear solvers and suggests new possibilities for designing iterative least squares solvers as iterative refinement. Specifically, sketched iterative refinement (or Iterative Hessian Sketching) can be interpreted as a preconditioned Jacobi iteration using the sketched matrix.  Assuming the meta-algorithm has a linear form $ALG^{meta}(A^\top b) = TA^\top b+q$, the sketched iterative refinement performs iteration $x_{i+1}  = (I-T^{-1}A^\top A)x_i+T^{-1}A^\top b$,
which is equivalent to Jacobi iteration with pre-conditer $T$. This indicates that the iterative refinement process implicitly acts as a preconditioning mechanism, enjoying the same convergence guarantees as described in \cite{rokhlin2008fast}. Moreover, this new understanding of iterative refinement allows for a more detailed analysis of numerical stability of the solver shown in Section \ref{section:SIRRfloat}.

\paragraph{Convergence of Iterative and Recursive refinement}
In this section we demonstrate the convergence of $\|x-x^\ast\|$.

\begin{theorem}[Convergence of Iterative/Recursive Refinement]
\label{convergence}
Suppose that the meta-algorithm has a linear form $ALG^{meta}(A^\top b) = TA^\top b+q$, then SIR and SRR are convergent if and only if $\rho(I-TA)<1$, with

\begin{itemize}
\setlength{\itemsep}{0pt}
\setlength{\parsep}{0pt}
\setlength{\parskip}{0pt}
\item {\small$\|\text{SIR}^{\text{meta}}_t(A^\top b)-x^\ast\|\leq \|\text{SIR}^{\text{meta}}_0(A^\top b)-x^\ast\|e^{-\alpha t}$},
\item {\small$\|\text{SRR}^{\text{meta}}_t(A^\top b)-x^\ast\|\leq \|\text{SRR}^{\text{meta}}_0(A^\top b)-x^\ast\|e^{-\alpha 2^t}$}

\end{itemize} 
where $\alpha = -ln (\rho (I-TA))$ and  $x^\ast$ is the true solution which satisfies $x^\ast=\arg\min_x \|Ax-b\|$.
\end{theorem}
\begin{remark}[Selection of Meta-Algorithm]\label{remark:select}
If one use the standard sketch-and-solve algorithm as the meta-algorithm, $t$-th iteration of SIR algorithm convergence at speed  $(\frac{1}{(1-\eta)^2}-1)^t$ for a sketching matrix with distortion $\eta$ where $\eta\in(0,1)$.  This means necessary sketching dimension depends on the intrinsic complexity of the problem. The algorithm would diverge if the “sufficient sketching dimension” condition
is violated \cite{pilanci2016iterative,wang2017sketching}. To remove such condition, we consider a 2-step Krylov-based sketch-and-solve solver as the meta-algorithm, now the $t$-th iteration of SIR algorithm convergence at speed $\min\{\eta^k,\frac{1}{\eta^k}\}$ which removes the requirement that $\eta<1$ (detailed proof shown in Appendix \ref{proofconvergence}).  We use the 2-step Krylov solver both for the stability analysis in Section \ref{section:SIRRfloat}  and the implementation in Section \ref{section:numerical}.
\end{remark}

    \section{Randomized Solver As Iterative Refinement}

In this section, we present a novel approach for constructing a fast and stable randomized least squares solver by iteratively refining an approximate solver which we call a meta-algorithm, \emph{e.g.} sketch-and-apply or early stopped iterative randomized solver. We introduce two ways to do the refinement: iterative refinement and recursive refinement. Both refinement process starts from a meta-algorithm and improve the previous solution by correcting it based on the residual error. The key difference between iterative and recursive refinement processes is that iterative refinement improves the solution by applying the meta-algorithm at each step to correct the residual, while recursive refinement refines the solution by repeatedly applying the same current solver to the residual error.

\begin{algorithm}

\caption{\textbf{S}ketched \textbf{I}terative \textbf{R}efinement}

\label{alg:Iterativerefinement}
    \begin{tcolorbox}
    \hrule\hrule
    \textbf{SIR}: Sketched Iterative Refinement
    \hrule
    \SetKwInOut{Input}{Input}
    \SetKwInOut{Output}{Output}
    \SetKw{KwBy}{by}
    \SetKw{KwReturn}{Return}
    \SetKw{KwVia}{Via}
    \Input{1}
    \Output{2}
    \hrule\hrule
    
    \If{$N=0$}{\KwReturn $\text{SIR}^{\text{meta}}_0(b)$ \KwVia a meta-algorithm $\text{SIR}^{\text{meta}}_0(b)=\text{ALG}^{\text{meta}} (A^\top b)$\;
    
    \Comment*[r]{Initialization via Meta-Algorithm}}
    \For{$i\gets1$ \KwTo $N$ \KwBy $1$}{
    $\text{SIR}^{\text{meta}}_i(b):=\text{SIR}^{\text{meta}}_{i-1}(b)+{\color{orange}\text{ALG}^{\text{meta}}} (A^\top (b-A \cdot \text{SIR}^{\text{meta}}_{i-1}(b)))$\Comment*[r]{Iterative Refinement via Meta-Algorithm}
    }
    \KwReturn $\text{SIR}^{\text{meta}}_N(b)$
    \end{tcolorbox}
    \end{algorithm}

\subsection{Iterative and Recursive refinement}
\label{section:IRandRR}

\paragraph{Iterative Refinement} Iterative refinement \cite{wilkinson2023rounding,moler1967iterative,wu2023stable} is the classical approach to improving the quality of a computed solution in numerical linear
algebra. The idea of iterative refinement is simple, to improve the quality of an approximate solution $x_i$, solve for the error $\delta x_i=x-x_i$ via approximately solving $\delta x_i:=\arg\min_{\delta x_i}\|b-Ax_i-A\delta x_i\|$. Classically, the inexact solve used in the refinement step is a classical direct solver such as QR factorization computed in lower numerical precision (i.e., single precision), and all the other steps are performed in higher precision (e.g., double precision) \cite{golub1966note,carson2024comparison}. In our paper, we design an iterative algorithm, where each step incorporates the concept of iterative refinement, using a fast randomized linear solver to approximately solve the system. The algorithm is detailed in Algorithm \ref{alg:Iterativerefinement}.

\paragraph{Recursive Refinement} We also introduce a novel way to implement an iterative refinement process which we call it (sketched) recursive refinement approach. Sketched Recursive Refinement process also iteratively refines the solution by incorporating corrections from previous iterations. Different from iterative refinement which updates the current solution by applying a fixed procedure to adjust the solution, recursive refinement refers back to itself to perform the next step and solve the problem in a nested fashion. The algorithm is detailed in Algorithm \ref{alg:recursiverefinement}. Later, we demonstrate that recursive refinement is simply a reorganization of the computational steps in iterative refinement but the two types of refinement enjoy very different numerical stability behavior.

\begin{algorithm}
\caption{\textbf{S}ketched \textbf{R}ecursive \textbf{R}efinement.}\label{alg:recursiverefinement}
\begin{tcolorbox}
\hrule\hrule
\textbf{SRR}: Sketched Recursive Refinement
\hrule
\SetKwInOut{Input}{Input}
\SetKwInOut{Output}{Output}
\SetKw{KwBy}{by}
\SetKw{KwReturn}{Return}
\SetKw{KwVia}{Via}
\Input{1}
\Output{2}
\hrule\hrule

\If{$N=0$}{\KwReturn $\text{SRR}_0(b)$ \KwVia meta-algorithm $\text{ALG}^{meta}(A^\top b)$\;}
\For{$i\gets1$ \KwTo $N$ \KwBy $1$}{
    $\text{SRR}_i(b):=\text{SRR}_{i-1}(b)+{\color{orange}\text{SRR}_{i-1}}(A^\top b-A^\top  A\cdot \text{SRR}_{i-1}(b))$\Comment*[r]{Recursive Refinement}
    }
\KwReturn $\text{SRR}_N(b)$
\end{tcolorbox}
\end{algorithm}

\paragraph{Recursive Refinement as Reorganizing Computation} We would like to point out that Recursive refinement and Iterative refinement perform the same if one uses exact arithmetic. With a linear meta-algorithm, \emph{i.e.} $ALG^{meta}(A^\top b)$ can be represented as $TA^\top b+q$ for some matrix $T$ which includes most useful randomized solver such as Sketch-and-Apply, the results of $\text{SIR}_{N}(b)$ and $\text{SRR}_{\log_2 N}(b)$ are the same and both can be presented in the same form as geometric series as $x=\sum_{i=0}^{N} (I-TA)^{i}Tb$ with same amount of compute $O(Nmn)$. This means that Recursive Refinement is just a reorganization of computation order in the Iterative Refinement procedure and would generate the same computational result if one use exact arithmetic.  However, in the following discussion, we show that Recursive Refinement and Iterative Refinement behave very differently when using a floating point arithmetic.

\paragraph{Equivalence between Iterative Refinement and Sketch-and-Precondition}
\label{section:equal}

Iterative Refinement (Iterative Hessian Sketching) and the Sketch-and-Precondition approach are commonly regarded as two distinct methodologies for designing iterative randomized least squares solvers. In this remark, we demonstrate the surprising equivalence between sketched iterative refinement and the sketch-and-precondition method. This insight provides a unified perspective on modern randomized linear solvers and suggests new possibilities for designing iterative least squares solvers as iterative refinement. Specifically, sketched iterative refinement (or Iterative Hessian Sketching) can be interpreted as a preconditioned Jacobi iteration using the sketched matrix.  Assuming the meta-algorithm has a linear form $ALG^{meta}(A^\top b) = TA^\top b+q$, the sketched iterative refinement performs iteration $x_{i+1}  = (I-T^{-1}A^\top A)x_i+T^{-1}A^\top b$,
which is equivalent to Jacobi iteration with pre-conditer $T$. This indicates that the iterative refinement process implicitly acts as a preconditioning mechanism, enjoying the same convergence guarantees as described in \cite{rokhlin2008fast}. Moreover, this new understanding of iterative refinement allows for a more detailed analysis of numerical stability of the solver shown in Section \ref{section:SIRRfloat}.

\paragraph{Convergence of Iterative and Recursive refinement}
In this section we demonstrate the convergence of $\|x-x^\ast\|$.

\begin{theorem}[Convergence of Iterative/Recursive Refinement]
\label{convergence}
Suppose that the meta-algorithm has a linear form $ALG^{meta}(A^\top b) = TA^\top b+q$, then SIR and SRR are convergent if and only if $\rho(I-TA)<1$, with

\begin{itemize}
\setlength{\itemsep}{0pt}
\setlength{\parsep}{0pt}
\setlength{\parskip}{0pt}
\item {\small$\|\text{SIR}^{\text{meta}}_t(A^\top b)-x^\ast\|\leq \|\text{SIR}^{\text{meta}}_0(A^\top b)-x^\ast\|e^{-\alpha t}$},
\item {\small$\|\text{SRR}^{\text{meta}}_t(A^\top b)-x^\ast\|\leq \|\text{SRR}^{\text{meta}}_0(A^\top b)-x^\ast\|e^{-\alpha 2^t}$}

\end{itemize} 
where $\alpha = -ln (\rho (I-TA))$ and  $x^\ast$ is the true solution which satisfies $x^\ast=\arg\min_x \|Ax-b\|$.
\end{theorem}

\begin{remark}[Selection of Meta-Algorithm]\label{remark:select}
If one use the standard sketch-and-solve algorithm as the meta-algorithm, $t$-th iteration of SIR algorithm convergence at speed  $(\frac{1}{(1-\eta)^2}-1)^t$ for a sketching matrix with distortion $\eta$ where $\eta\in(0,1)$.  This means necessary sketching dimension depends on the intrinsic complexity of the problem. The algorithm would diverge if the “sufficient sketching dimension” condition
is violated \cite{pilanci2016iterative,wang2017sketching}. To remove such condition, we consider a 2-step Krylov-based sketch-and-solve solver as the meta-algorithm, now the $t$-th iteration of SIR algorithm convergence at speed $\min\{\eta^k,\frac{1}{\eta^k}\}$ which removes the requirement that $\eta<1$ (detailed proof shown in Appendix \ref{proofconvergence}).  We use the 2-step Krylov solver both for the stability analysis in Section \ref{section:SIRRfloat}  and the implementation in Section \ref{section:numerical}.
\end{remark}

\section{Fast and Stable Solver via Iterative and Recursive refinement}

To construct a fast and stable randomized solver, we use Sketched Recursive Refinement as the meta-algorithm for a Sketched Iterative Refinement process. We call our algorithm Sketched Iterative and Recursive Refinement (SIRR) {which is shown as algorithm \ref{alg:SIRR} in the appendix.} We also theoretically show that both iterative and recursive refinement are essential to achieve backward stability. The theoretical finding is also verified numerically in Section \ref{section:numerical}.

\subsection{Sketched Iterative and Recursive Refinement}

\paragraph{SIRR is Fast} In this section, we first show that SIRR converges fast with a computational complexity at \alertinline{$O(n^3+mn)$}. Note that SIRR is a composite of meta-algorithm, so we examine the computational complexity and average convergence rate of meta-algorithm to show the whole computational complexity of SIRR. 

Suppose that $A\in \mathbb{R}^{m\times n},b\in \mathbb{R}^{m\times 1}$ and upper triangular matrix $R\in \mathbb{R}^{n\times n}$, the computation of computational complexity follows:

\begin{itemize}
\setlength{\itemsep}{0pt}
\setlength{\parsep}{0pt}
\setlength{\parskip}{0pt}
    \item matrix and vector multiplication $A^\top b$: $O(mn)$
    \item solving triangular system $R^{-1}z$ for $z\in \mathbb{R}^{n\times 1}$: $O(n^2)$
    \item conducting QR factorization of $SA$: $O(sn^2)=O(n^3log(n))$
\end{itemize}

For sketch-and-solve meta solver, the computational complexity is $O(n^2+mn)$ and the convergence rate is $\frac{1}{(1-\eta)^2}-1$. To reach machine precision, the iteration step is at most $O(log(\frac{1}{u}))$, thus the total computational complexity of SIR, SRR and SIRR are all $O(n^3log(n)+log(\frac{1}{u})(n^2+mn))$ matches the fast randomzied least square solvers such as Blendenpik \cite{avron2010blendenpik} and FOSSILS \cite{epperly2024fast}.

\subsection{SIRR in Floating Point Arithmetic}

\label{section:SIRRfloat}
As shown in Section \ref{section:IRandRR}, SIRR is the same as SIR in exact arithmetic. In this section, we study the stability results for SIR, SRR and SIRR when one implements them in floating point arithmetic.

\subsubsection{SIRR is forward Stable}

In this section, we first prove that SIRR solver is forward stable, \emph{i.e.} both the forward error $\|\hat x - x^\ast \|$ and the residual error $\|A\hat x - Ax^\ast \|$ converge geometrically for SIRR implemented in floating point arithmetic.

\begin{definition}[Forward Stability] 
A least-squares solver is forward stable if the computed solution $\hat{x}$ satisfies 

    \begin{align*}
        \|\hat{x}-x^*\|\leq \epsilon(\kappa\|x^*\|+\frac{\kappa^2}{\|A\|}\|r^*\|),
   \end{align*}

and is strongly forward stable if $\hat{x}$ satisfies
\begin{align*}
    \|A(\hat{x}-x^*)\|\leq \epsilon(\kappa\|r^*\|+\|A\|\|x^*\|),
    \end{align*}
    where $x^*$ is the exact solution, $r^*=b-Ax^*$ and $\epsilon\lesssim n^\frac{3}{2}$
\end{definition}

\begin{remark} This is the best error one can expect to achieve due to  Wedin’s theorem \cite{wedin1973perturbation}, where one always solving a perturbed problem $\argmin_{y\in \mathbb{R}^{n}}\|(b+\delta b)-(A+\delta A)y\|$  in floating point arithmetic, where $\|\delta A\|\leq \epsilon\|A\|$, $\|\delta b\|\leq\epsilon \|b\|$. 
\end{remark}

\begin{theorem}
\label{forwardstable}
For SIRR with meta-algorithm $\text{ALG}^{meta}(\cdot)$, which solves problem $ x=\arg\min_y \|(A^\top A)y-r_A\|$, satisfying $\text{ALG}^{meta}(r_A) = (A^\top A)^{-1}r_A+c\|\hat{R}^{-\top}r_A\|\hat{R}^{-1}e$ where $\|e\|\lesssim 1$ and $c<1$, the result $\hat{x}$ of SIRR is strongly forward stable, which satisfies

\begin{align*}
    \|\hat{x}-x^*\|\lesssim n^{\frac{3}{2}}(u\kappa\|x^*\|+\frac{u\kappa^2}{\|A\|}\|r^*\|),\\
    \|A(\hat{x}-x^*)\|\lesssim n^{\frac{3}{2}}(u\kappa\|r^*\|+u\|A\|\|x^*\|),
\end{align*}

\end{theorem}
With strongly forward stable, we can expect a non-pathological rounding error $\|\hat{x}\|\geq \|x^*\|+n^{\frac{3}{2}}\frac{u\kappa^2}{\|A\|}\|r^*\|.$

\subsubsection{SIRR is Backward Stable}

In this section, we provide the theoretical analysis (following $\alpha-\beta$ accuracy framework from \cite{epperly2024fast}) showing that the Sketched Iterative and Recursive Refinement (SIRR) is provable backward stable when implemented in floating point arithmetic. To do this, we first find the requirement that the meta-solver of the sketched iterative refinement needs to satisfy that can make SIR solver backward stable. Then we prove that Sketched Recursive Refinement can provably meet these requirements. 

\begin{theorem}
    \label{strong solver}
    For simplicity, denote $\max\{u\kappa,\frac{1}{\kappa n^{\frac{3}{2}}}\}$ as $\tilde{\kappa}^{-1}$. Suppose that $un^{\frac{3}{2}}\|x^\ast\|\le\|b-Ax^\ast\|$ and the single step meta-solver $ALG(z)$ is $n^{\frac{3}{2}}\tilde{\kappa}^{-1}(\|Ax_z^*\|+u\kappa\|z\|)-un^{\frac{3}{2}}(\|Ax_z^*\|+\|z\|)$ accurate where $x_z^\ast$ is the true solution of the least square problem, \emph{i.e.} $\argmin_x \|Ax-z\|$. Then SIR solver $x_{i+1}=ALG(b-Ax_i)$ will converge to a $(un^{3}\kappa\tilde{\kappa}^{-1}\|b-Ax^\ast_b\|+un^{\frac{3}{2}}\|x^\ast\|)-(un^{\frac{3}{2}}\|b-Ax^\ast_b\|+u^2n^3\|x^\ast\|)$ accurate solution
    which indicate a \textbf{backward stable} result by Lemma \ref{becondition}.
\end{theorem}
\begin{remark}
\label{nremark} Since SIR/SIRR solver enjoys non-pathological rounding error assumption $\|x^\ast\|+n^{\frac{3}{2}}\kappa^2u\|b-Ax^\ast\|\lesssim\|\hat{x}\|$ (Theorem \ref{forwardstable}), we have ${\scriptsize\overbrace{(n^{\frac{3}{2}}u+n^3\kappa^2u^2)}^{{un^{3}\kappa\tilde{\kappa}^{-1}}}}\|b-Ax^\ast\|+un^{\frac{3}{2}}\|x^\ast\|
        \lesssim un^{\frac{3}{2}}\|b-Ax^\ast\|+un^{\frac{3}{2}}\|\hat{x}\|\lesssim un^{\frac{3}{2}}\|b-A\hat{x}\|+un^{\frac{3}{2}}\|\hat{x}\|
        \lesssim un^{\frac{3}{2}}+un^{\frac{3}{2}}\|\hat{x}\|$
    and $un^{\frac{3}{2}}\|b-Ax^\ast\|+u^2n^3\|x^\ast\|\lesssim un^{\frac{3}{2}}\|b-Ax^\ast\|\lesssim un^{\frac{3}{2}}\|b-A\hat{x}\|$ based on the assumption that $un^{\frac{3}{2}}\|x^\ast\|\lesssim\|b-Ax^\ast\|$. By lemma \ref{becondition}, the solution has backward error $Be(\hat{x})\lesssim n^2u$ which indicates a backward stability result and \emph{aligns the backward error estimation for QR-based solver \cite[Theorem 19.5]{higham2002accuracy} which also dependency of matrix size at $n^2$}.
\end{remark}

Then we study the stability result of SRR implemented in floating point arithmetic. We show that SRR can be backward stable only when $\frac{\|b-Ax^\ast\|}{\|Ax^\ast\|} = O(1)$ and is not backward stable when the residual $\|b-Ax^\ast\|$ is small.  However, SRR provides an approximate solver that satisfies the assumption we require for the meta-algorithm in the backward stable result in Theorem \ref{strong solver}.

\begin{theorem}
    \label{iterativeIDS}
    {\color{red} }For meta-algorithm $SRR_0(\cdot)$, which solves problem $ x=\arg\min_y \|(A^\top A)y-r_A\|$, satisfying $SRR_0(r_A) = x^*+(a_1\|x^*\|+a_2\|Ax^*\|)\hat{R}^{-1}e_1+(b_1\|x^*\|+b_2\|Ax^*\|)(A^\top A)^{-1}e_2$, where $x^* = (A^\top A)^{-1}r_A$ and 
    \begin{align*}
        \kappa a_1+a_2\asymp c, \kappa^2 b_1+\kappa b_2\asymp c, \|e_{1,2}\|\lesssim 1,
    \end{align*}
    and $N=O(log_2(\frac{log(\tilde{\kappa}^{-1}n^{\frac{3}{2}})}{log(c)}))$, the solution of corresponding $\text{SRR}_N$ is $O(n^{\frac{3}{2}}(\tilde{\kappa}^{-1}\|Ax^\ast\|+u\|x^\ast\|+u\kappa\tilde{\kappa}^{-1}\|b\|),n^{\frac{3}{2}}(u\|Ax^\ast\|+u^2n^{\frac{3}{2}}\|x^\ast\|+u\|b\|))$-accurate. As $N\to \infty$, $\text{SRR}_N$ converges to a $(a,b)-$accurate solution $\text{SRR}_\infty(b)$ with
    \begin{align*}
        \hat{a}&\lesssim (u^2\kappa^2n^3\|Ax^\ast\|+un^{\frac{3}{2}}\|x^\ast\|+u^2\kappa^2n^3\|b\|),\\
        \hat{b}&\lesssim (un^{\frac{3}{2}}\|Ax^\ast\|+u^2n^3\|x^\ast\|+un^{\frac{3}{2}}\|b\|).
    \end{align*}
\end{theorem}
\begin{remark} Theorem \ref{iterativeIDS} indicates that SRR has the same backward error level as SIRR when $\frac{\|b-Ax^\ast\|}{\|Ax^\ast\|} = O(1)$. We verified numerically that SRR solver only is not backward stable when $\|b-Ax^\ast\|$ is large. The result is presented in Figure \ref{fig:errorsize}. This illustrates that our theoretical result for SRR is tight.
\end{remark}

Although $\text{SRR}_N(\cdot)$ is not backward stable on its own, it satisfies the requirements (for $u\|x^*\|\lesssim \tilde{\kappa}^{-1}\|Ax^*\|$ and $un^{\frac{3}{2}}\|x^*\|\lesssim \|Ax^*\|$) of the meta-algorithm in Theorem \ref{strong solver} to achieve a backward-stable SIR solver, as demonstrated in Theorem \ref{iterativeIDS}. This implies that by using SRR as the meta-algorithm for the SIR solver—\emph{i.e.}, the SIRR solver—it can be proven to be backward-stable, provided the meta-algorithm satisfies the conditions outlined in Theorem \ref{iterativeIDS}. Finally, we show that the two-step Krylov-based meta-algorithm, described in Remark \ref{remark:select}, meets the meta-algorithm criteria specified in Theorem \ref{iterativeIDS}.
\begin{lemma}
\label{metasolver}
    The result of 2-step Krylov-based meta-algorithm (Appendix Algorithm \ref{alg:krylov}) for solving $ x=\arg\min_y \|(A^\top A)y-r_A\|$ satisfies 
    \begin{align*}
        \hat x &= x^*+u\kappa n^{\frac{3}{2}}\|Ax^*\|\hat{R}^{-1}e_1+ un^\frac{3}{2}\|Ax^*\| (A^\top A)^{-1}e_2,
    \end{align*}
    where $x^*=(A^\top A)^{-1}r_A$ and $\|e_{1,2}\|\lesssim 1$. As a result, SIRR with Krylov-based meta-algorithm is backward stable.
\end{lemma}

\section{Numerical Experiments}
\label{section:numerical}

In this section, we compare SIR, SRR and SIRR solver to verify our theoretical findings. We also compare it with QR-based direct solver (\textit{mldivide} (MATLAB)) and FOSSILS in concurrent work \cite{epperly2024fast} to show that SIRR solver can beat the state-of-the-art randomized/direct solvers in realistic applications.

\paragraph{Error metrics} Following \cite{epperly2023fast,epperly2024fast}, we test three useful error metrics for all randomized least square solvers:

\begin{enumerate}
\setlength{\itemsep}{0pt}
\setlength{\parsep}{0pt}
\setlength{\parskip}{0pt}
    \item \textbf{Forward error.} The forward error quantifies how close the computed solution $\hat{x}$ is to the true solution $x$, \emph{i.e.} $\text{FE}(\hat{x}) := \frac{\|x - \hat{x}\|}{\|x\|}.$

    \item \textbf{Residual error.}  The (relative) residual error measures the \textit{suboptimality} of $\hat{x}$ as a solution to the least-squares minimization problem, \emph{i.e.} $\text{RE}(\hat{x}) := \frac{\|r(x) - r(\hat{x})\|}{\|r(x)\|}$.
  
    \item \textbf{Backward error.} The (relative) backward error \cite[Section 20.7]{higham2002accuracy} is $\text{BE}_b(\hat{x}) := \min_v \frac{\|\Delta A\|_F}{\|A\|_F}$ 
   where $\hat{x} = \arg\min_v \|b - (A + \Delta A)v\|$. 
\end{enumerate}

\paragraph{Experiment Setup} We adopt a similar setup to \cite{meier2024sketch,epperly2023fast} in most of experiments. We set $A\in R^{m\times n}$, sketching matrix $S\in R^{s\times m}$, and choose parameters \( \kappa \geq 1 \) for the condition number of \( A \) and \( \beta \geq 0 \) for the residual norm \( \|r(x)\| \). To generate \( A \), \( x \), and \( b \), do the following:

\begin{itemize}
\setlength{\itemsep}{0pt}
\setlength{\parsep}{0pt}
\setlength{\parskip}{0pt}    
\item Choose Haar random orthogonal matrices \( U = [U_1 \ U_2] \) in \( \mathbb{R}^{m \times m} \) and \( V \) in \( \mathbb{R}^{n \times n} \), and partition \( U \) so that \( U_1 \in \mathbb{R}^{m \times n} \).
    \item Set \( A := U_1 \Sigma V^T \) where \( \Sigma \) is a diagonal matrix with logarithmically equispaced entries between 1 and \( \frac{1}{\kappa} \).
    \item Form vectors \( w \) in \( \mathbb{R}^n \), \( z \) in \( \mathbb{R}^{m-n} \) with independent standard Gaussian entries.
    \item Define the solution \( x := \frac{w}{\|w\|} \), residual \( r(x) = \beta \cdot U_2 z / \|U_2 z\| \), and right-hand side \( b := Ax + r(x) \).
\end{itemize}

 We also experiment on kernel regression task, where we consider least-squares problems for fitting the SUSY dataset using a linear combination of kernel functions. Similar to \cite{epperly2023fast,epperly2024fast}, we generate real-valued least-squares problems of dimension $m = 10^6$ and $n \in [10^1, 10^3]$. 

\begin{figure}
    \centering
  {\scriptsize$$\|r(x)\|=10^{-1}$$}
\includegraphics[width=1in]{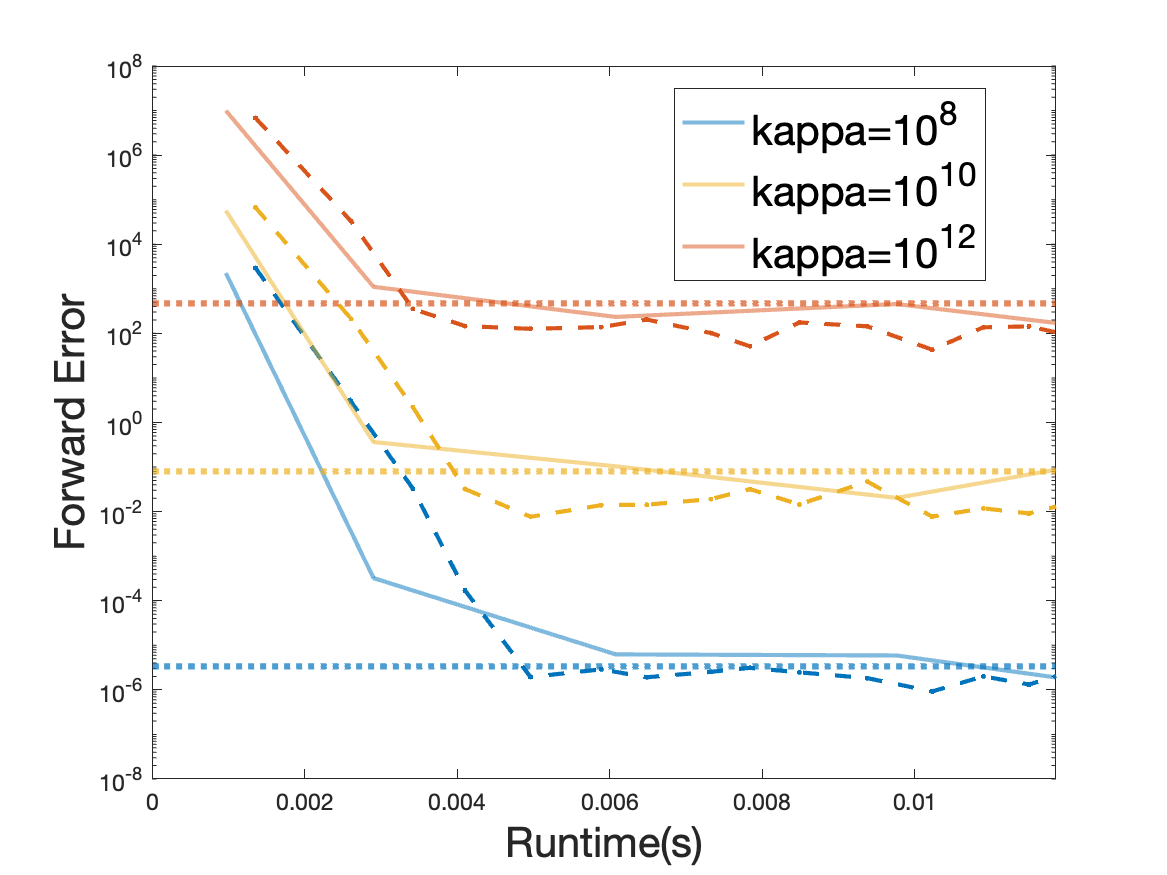}
    \includegraphics[width=1in]{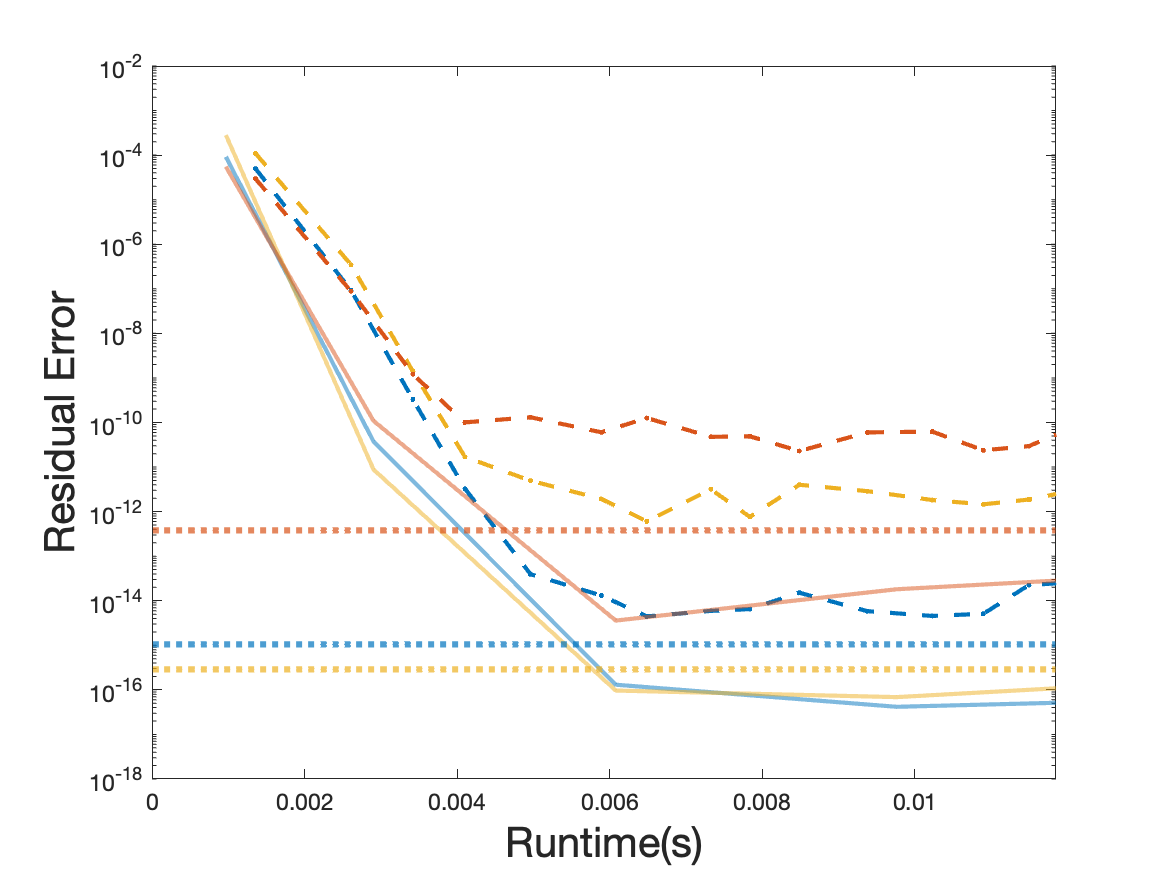}
    \includegraphics[width=1in]{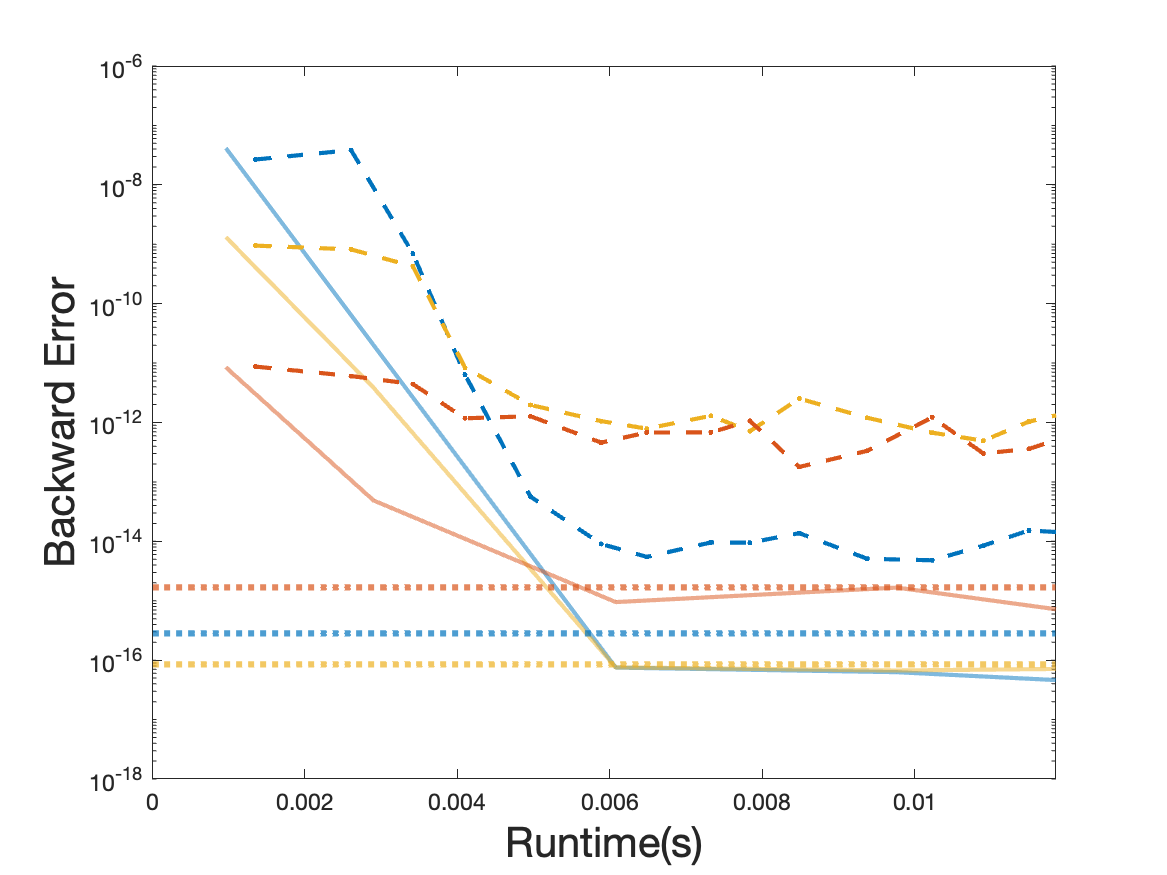}

 {\scriptsize$$\|r(x)\|=10^{-3}$$}

\includegraphics[width=1in]{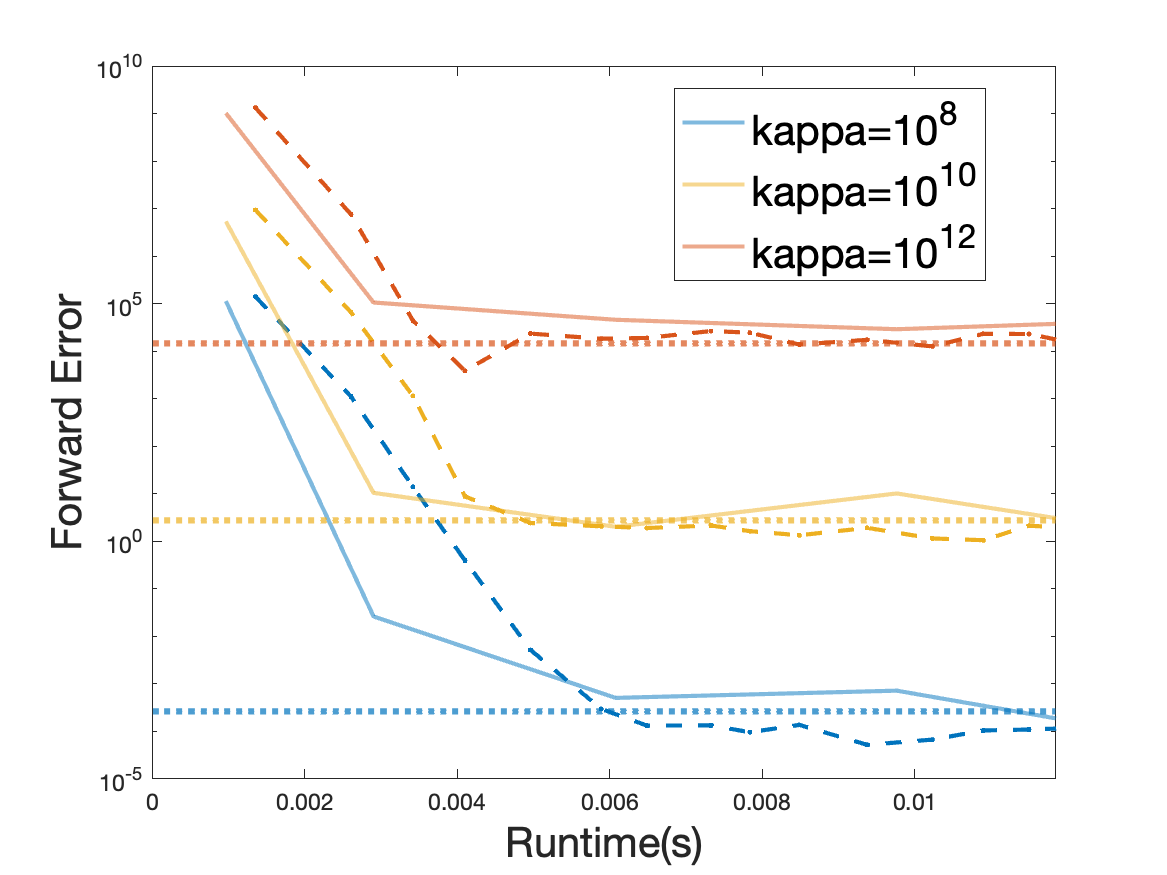}
    \includegraphics[width=1in]{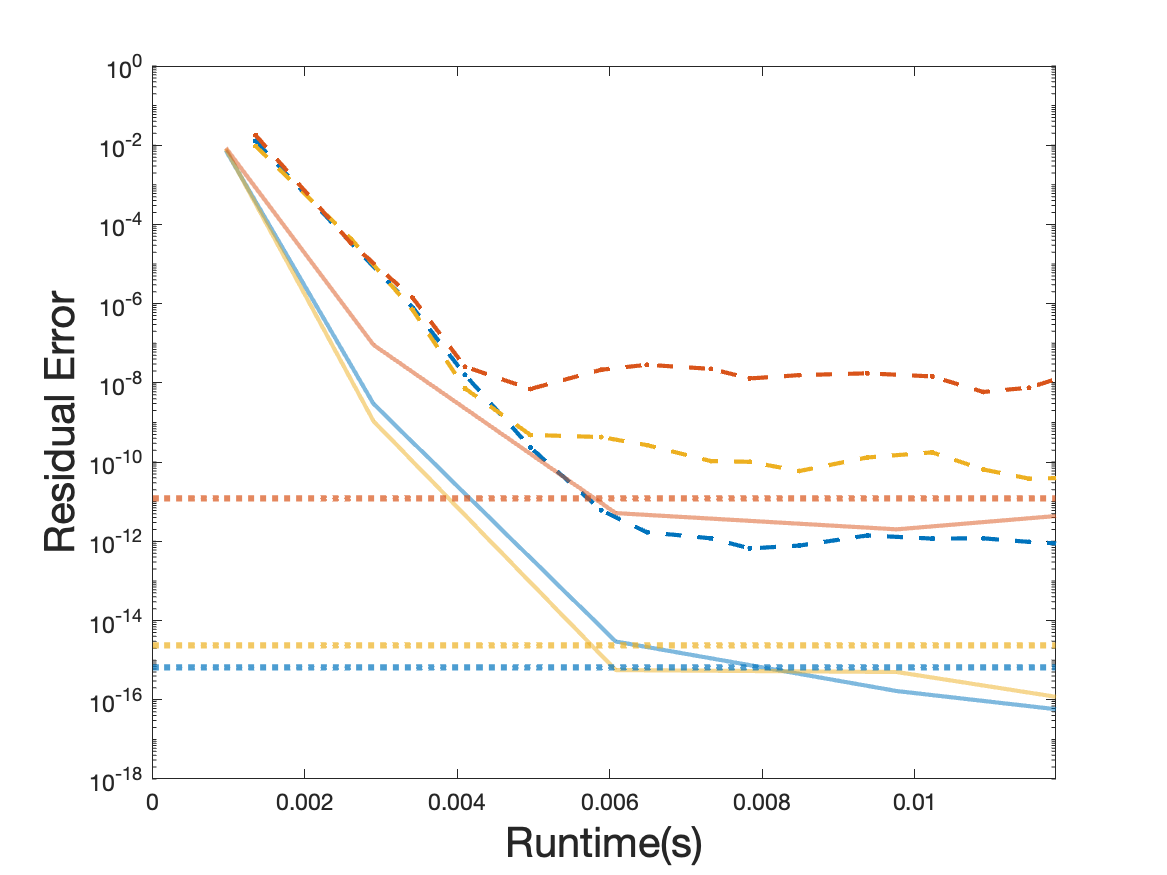}
    \includegraphics[width=1in]{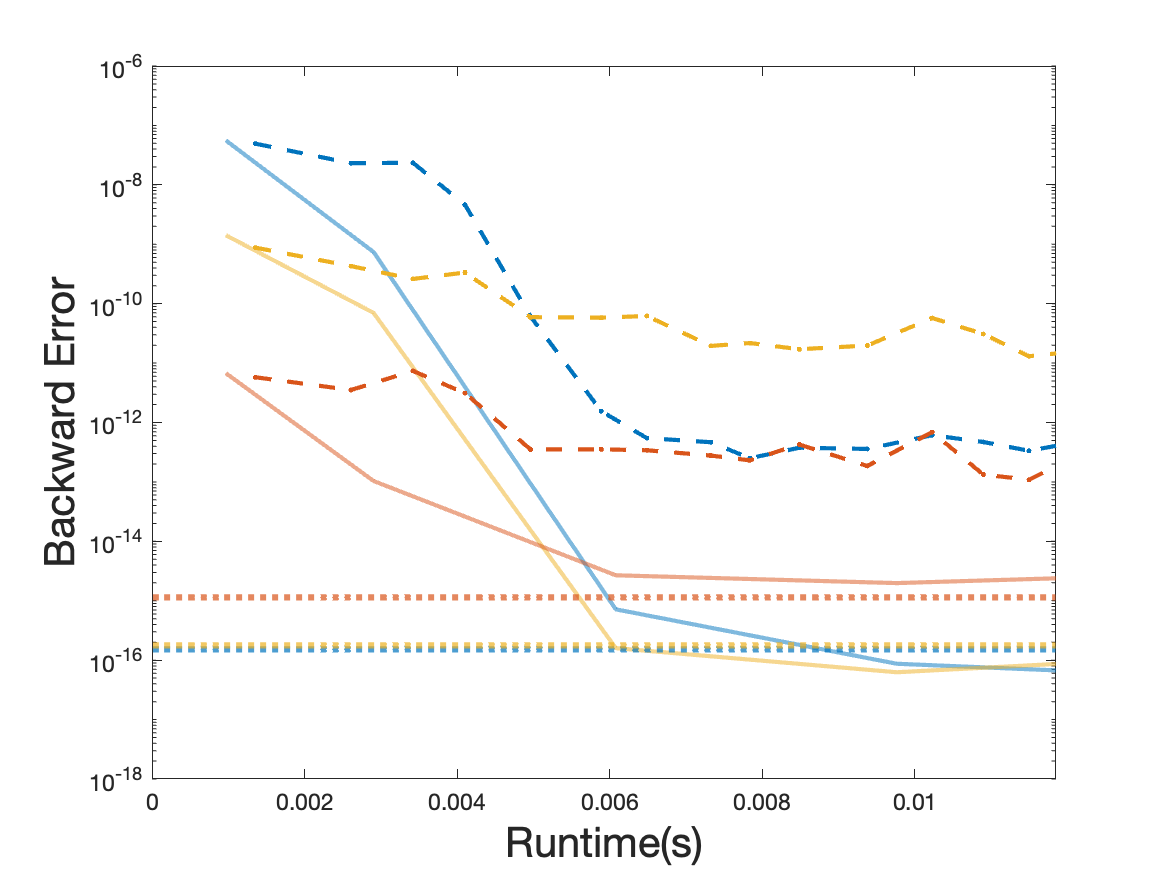}
    \caption{Results of SIRR with sketch and solve Initialization are shown as solid curve lines, with reference accuracy for MATLAB function \textit{A$\backslash$b} shown as dotted constant lines and IHS-Krylov shown as dotted curve lines }
\label{fig:idsreszeroinit}
\end{figure}

\begin{figure}
    \centering
    \includegraphics[width=1.5in]{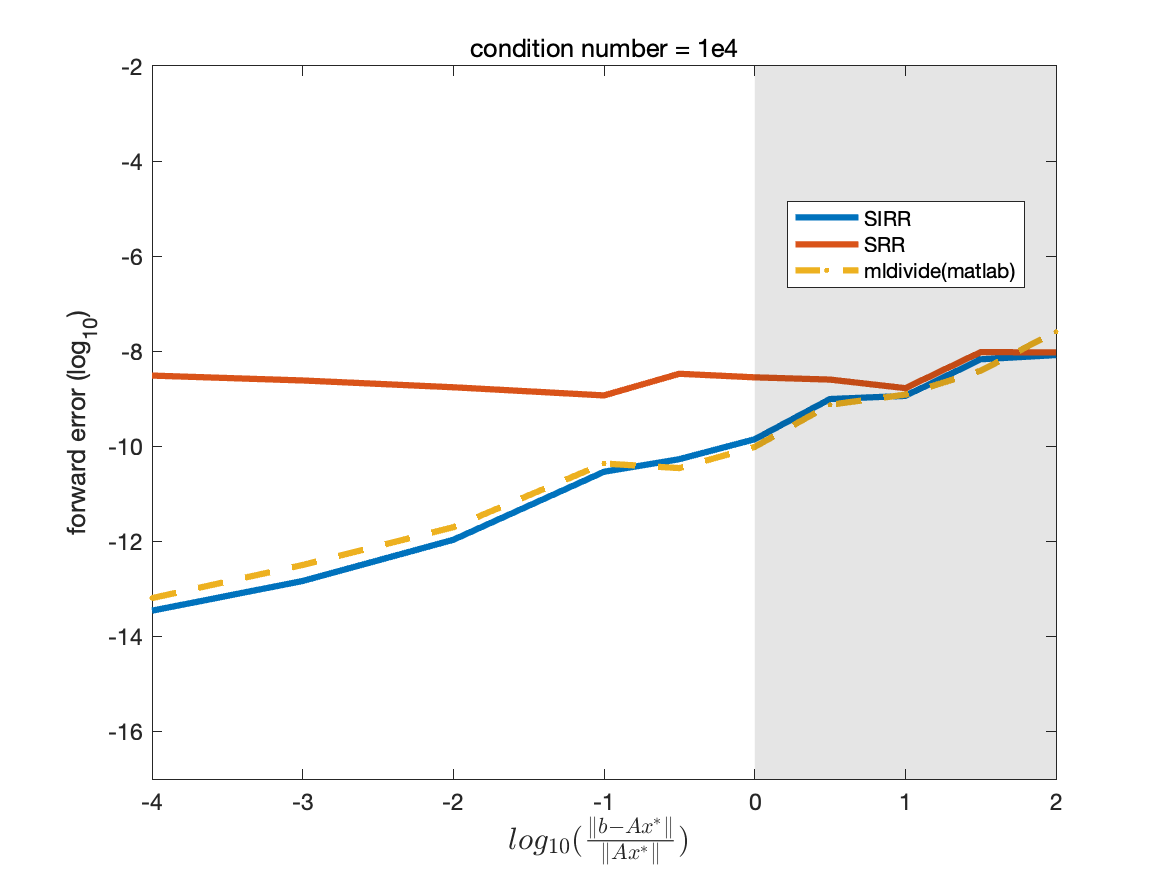}
    \includegraphics[width=1.5in]{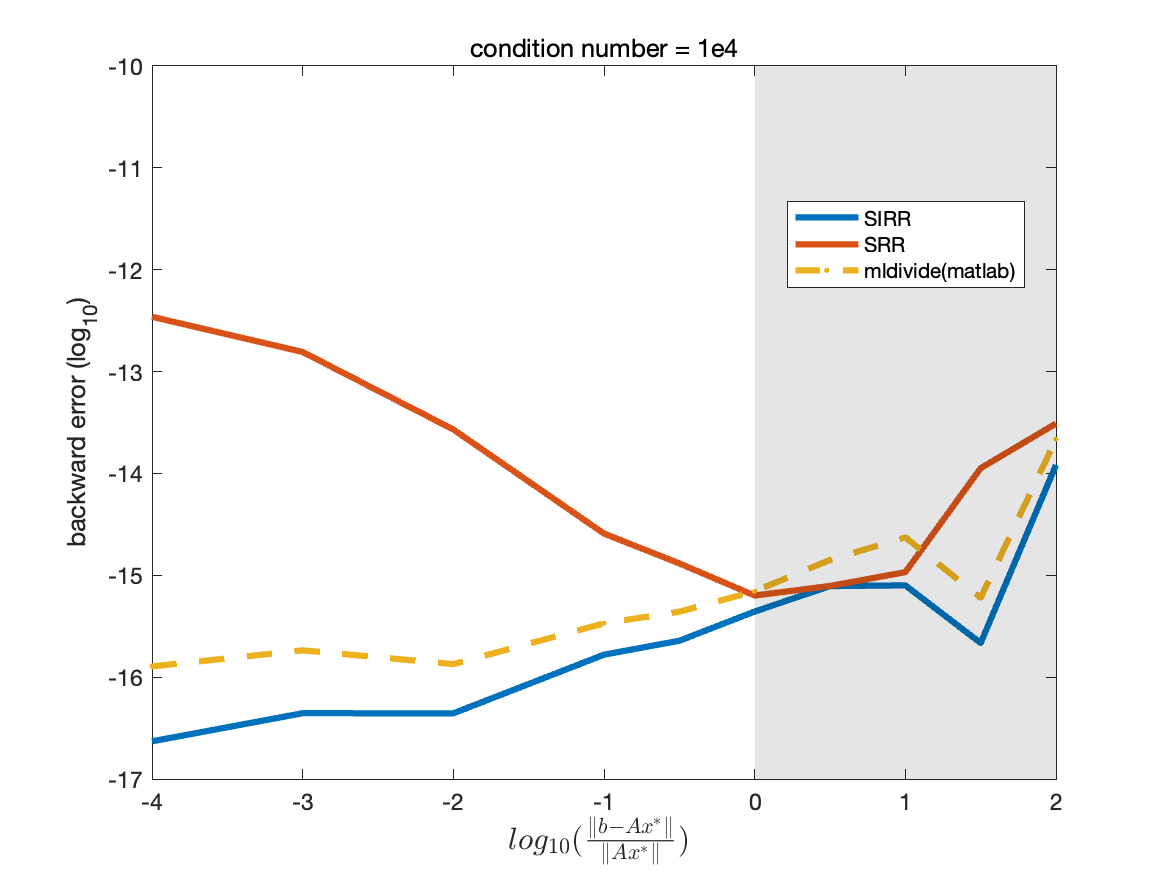}
    \includegraphics[width=1.5in]{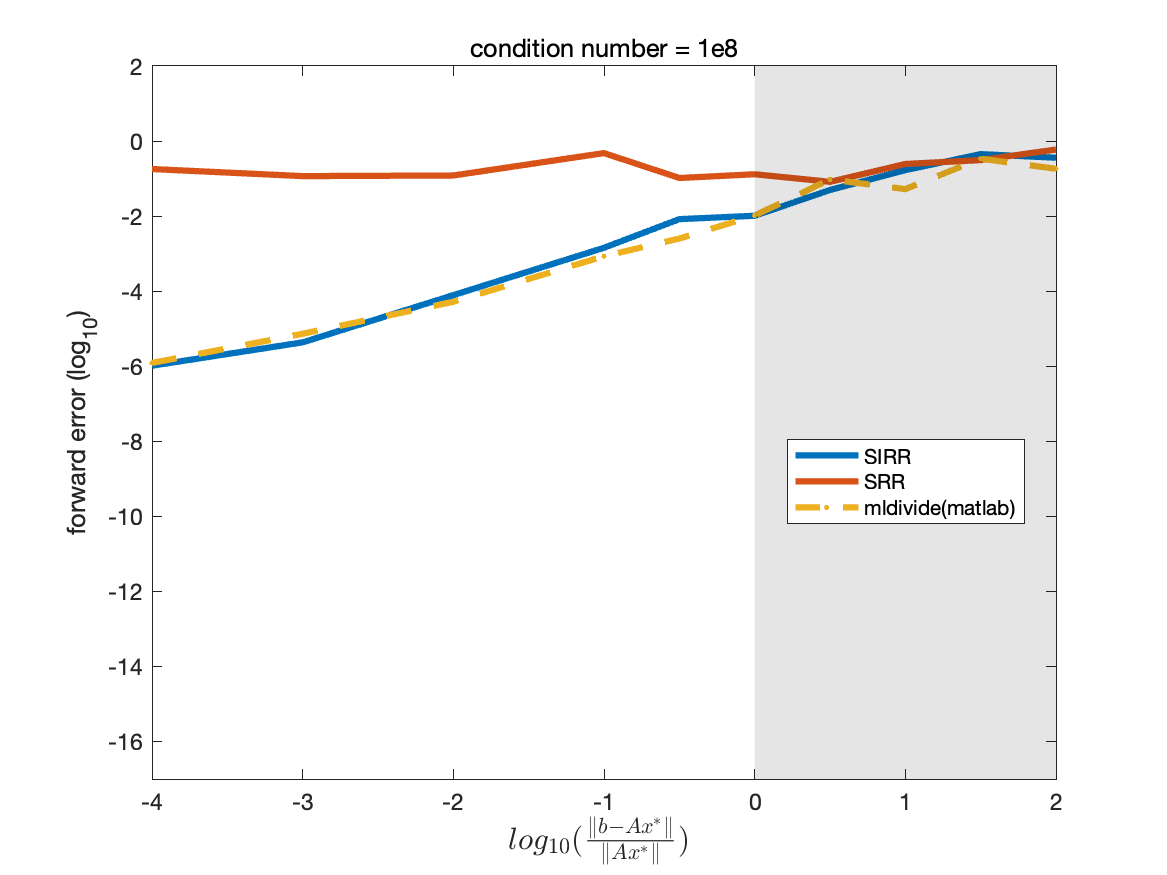}
    \includegraphics[width=1.5in]{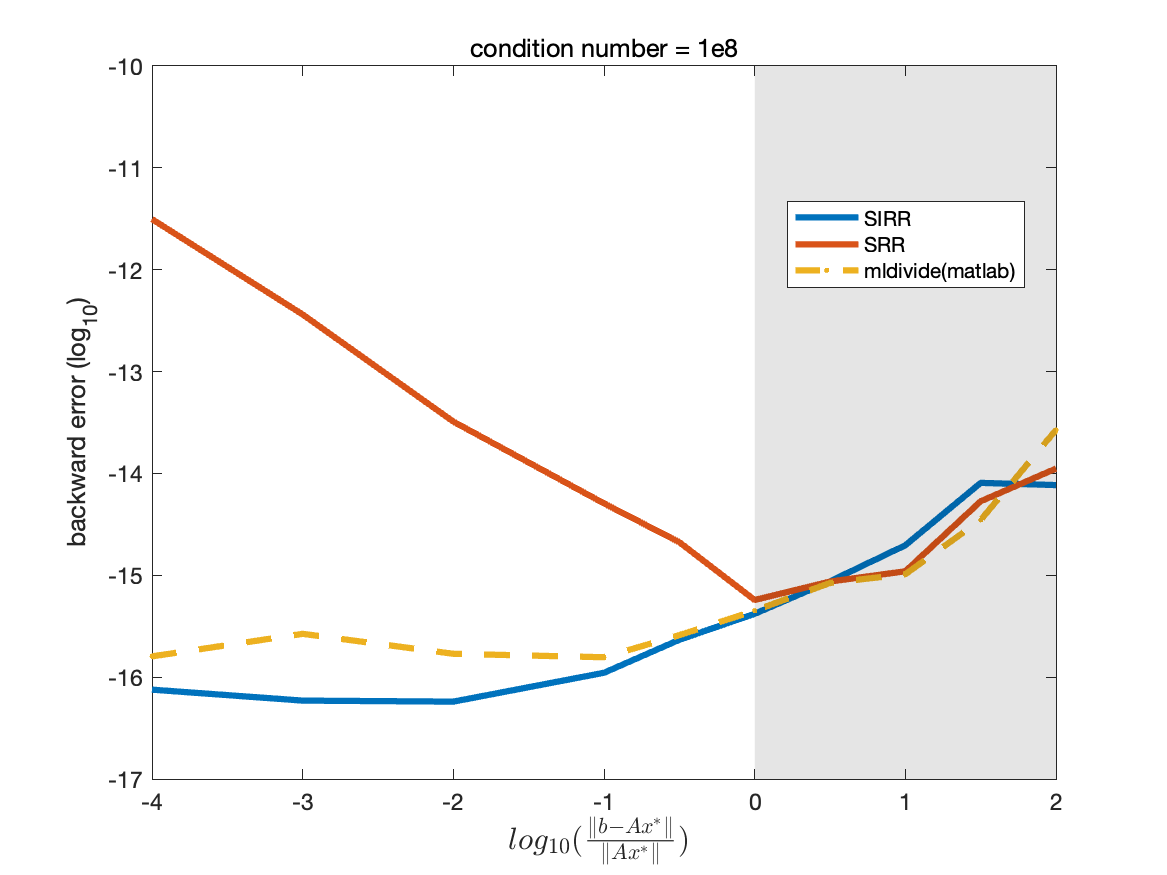}
    \includegraphics[width=1.5in]{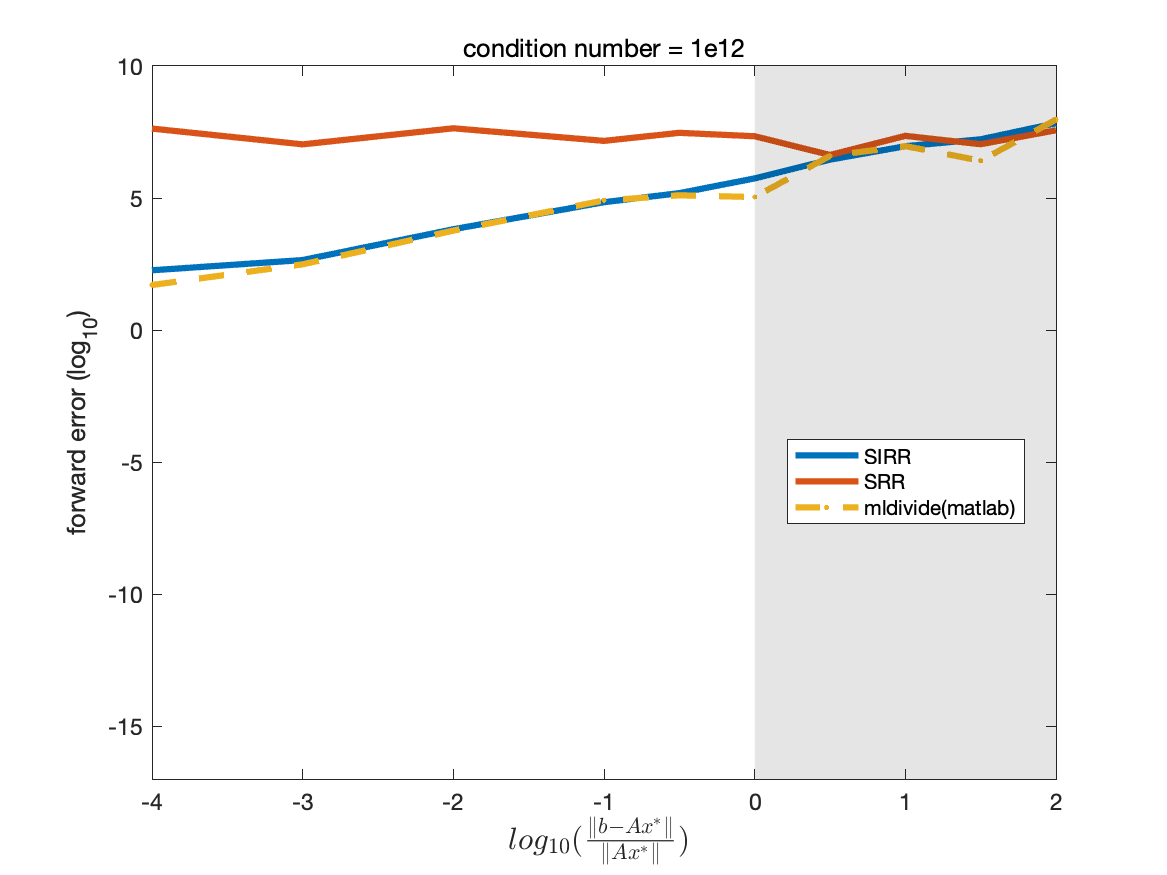}
    \includegraphics[width=1.5in]{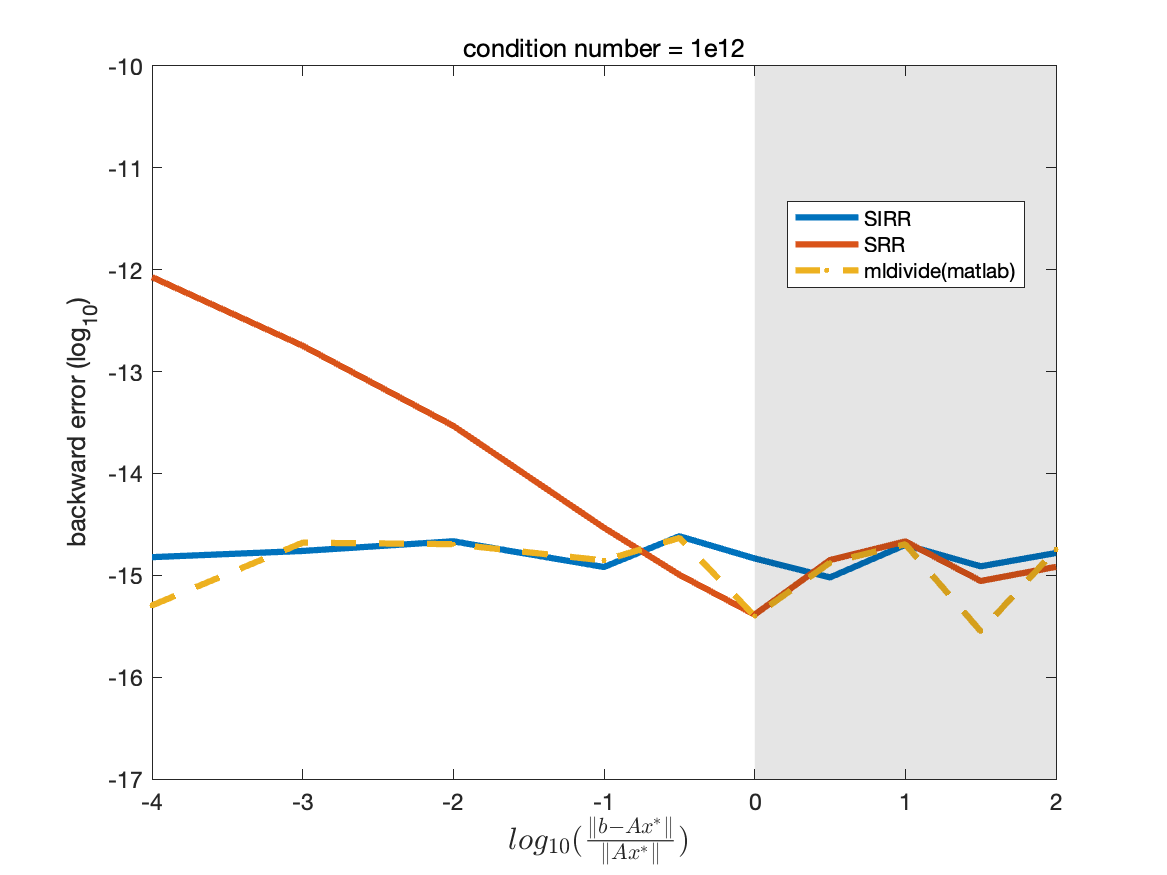}
    \caption{Forward error (left) and backward error (right) under different  $\|b-Ax^\ast\|/\|Ax^\ast\|$. SRR is not backward stable when $\|b-Ax^\ast\|$ is small while SIRR can achieve backward stable estimates for all cases. We also plotted the result for \textit{mldivide}(MATLAB) solver here for reference.}
\label{fig:errorsize}
\end{figure}

\paragraph{Both Iterative and Recursive Refinement is Essential} In this section, we conduct numerical experiments to demonstrate that both iterative and recursive refinement are essential for constructing a backward-stable solver. To illustrate this, we compare SIR, SRR, and SIRR, each using a two-step Krylov solver {as algorithm \ref{alg:krylov} in appendix.} as the meta-solver, under varying levels of condition numbers and $\frac{\|b - Ax^\ast\|}{\|Ax^\ast\|}$ to validate theorem \ref{iterativeIDS}. Figure \ref{fig:idsreszeroinit} shows that the SIR solver is not backward stable, while the SIRR solver achieves near machine-precision backward error. In a second experiment, we compare SIRR and SRR across different levels of residual size. Our theoretical results in theorem \ref{iterativeIDS} indicate that when the magnitude of the residual $\|b - Ax^\ast\|$ exceeds the signal $\|Ax^\ast\|$, SRR achieves the same backward stability as SIRR. However, SRR cannot achieve the same level of backward stability as SIRR when the residual $\|b - Ax^\ast\|$ is small. Figure \ref{fig:errorsize} confirms this result, showing that the backward error of SRR converges to that of SIRR in the white region and reaches the same level as SIRR in the grey region. In all experiments, we set $m=2000, n=50, s=200$.

\begin{figure}
    \centering
    \includegraphics[width=1.5in]{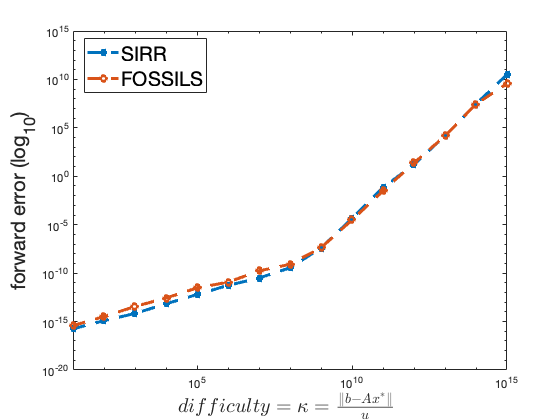}
    \includegraphics[width=1.5in]{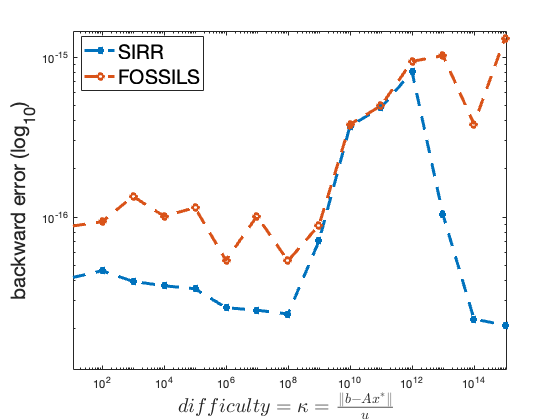}
  \caption{Comparing the Forward error (left) and backward error (right) of SIRR and FOSSILS on problems with different difficulties. SIRR has better backward stability in most situations and similar forward stability compared to FOSSILS.}
    \label{fig:SIRRvsFOSSILS}
\end{figure}

\paragraph{SIRR VS FOSSILS} We also compare our SIRR solver with FOSSILS in concurrent work \cite{epperly2024fast}\footnote{We use the code from \url{https://github.com/eepperly/Stable-Randomized-Least-Squares} for the FOSSILS algorithm. After an email exchange with authors of \cite{epperly2024fast}, we noticed that they have remarked that remark that, for smaller values of $d$,  a slightly larger value $\eta = 1.2 \sqrt(n/d)$ may fix the instability. However, this modification is not included in the GitHub version of the FOSSILS code. All results reported in our manuscript are based on the GitHub version of the solver.} by two experiments. In the first experiments, we adopt the same setting as \cite{epperly2024fast}, where a family of problems is generated of increasing difficulty, with condition number $\kappa$ and error size $\|b-Ax^\ast\|$ satisfying 
{\footnotesize$$\text{difficulty}=\kappa=\frac{\|b-Ax^\ast\|}{u}\in[10^0,10^{16}].
$$} We set $m=5000,n=200,s=600$ for problem size.Figure \ref{fig:SIRRvsFOSSILS} shows the forward error and backward error of SIRR and FOSSILS in problems of different difficulties, where both sketching algorithms have a similar forward stability while SIRR exhibits a better backward error performance.

\begin{figure*}
\label{exp:runtime}
    \centering

    \includegraphics[width=2in]{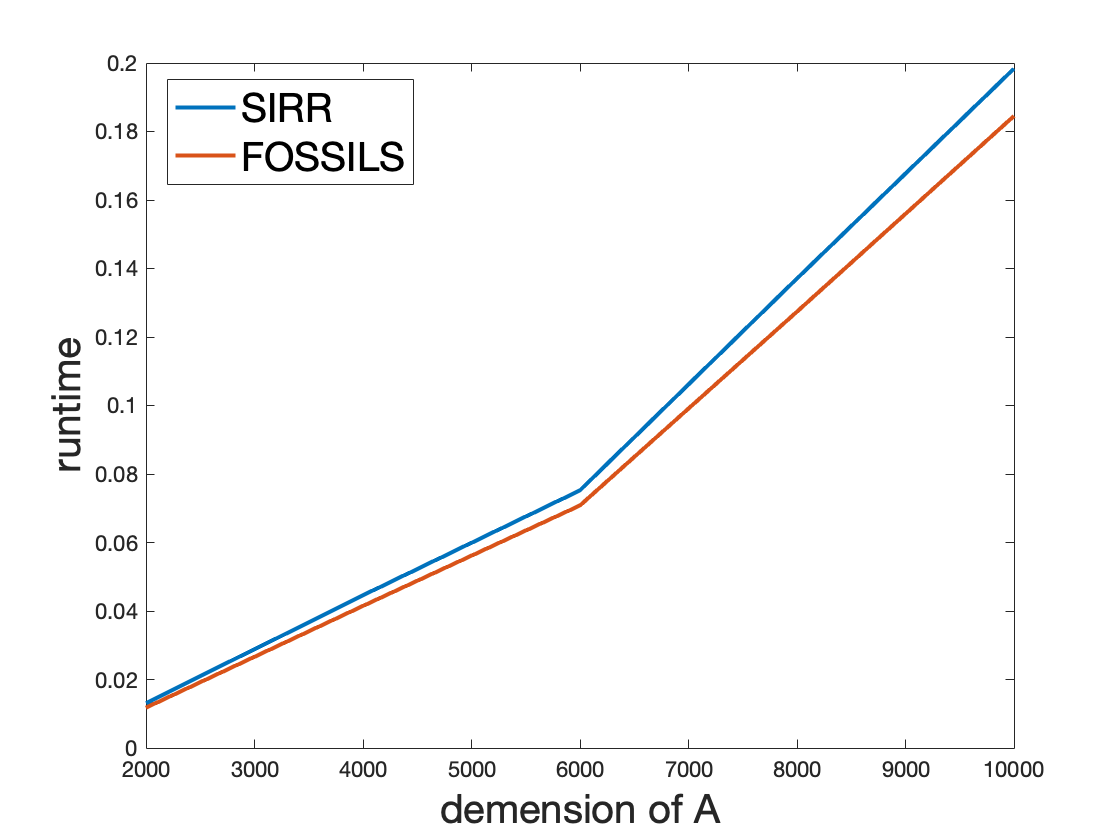}
    \includegraphics[width=2in]{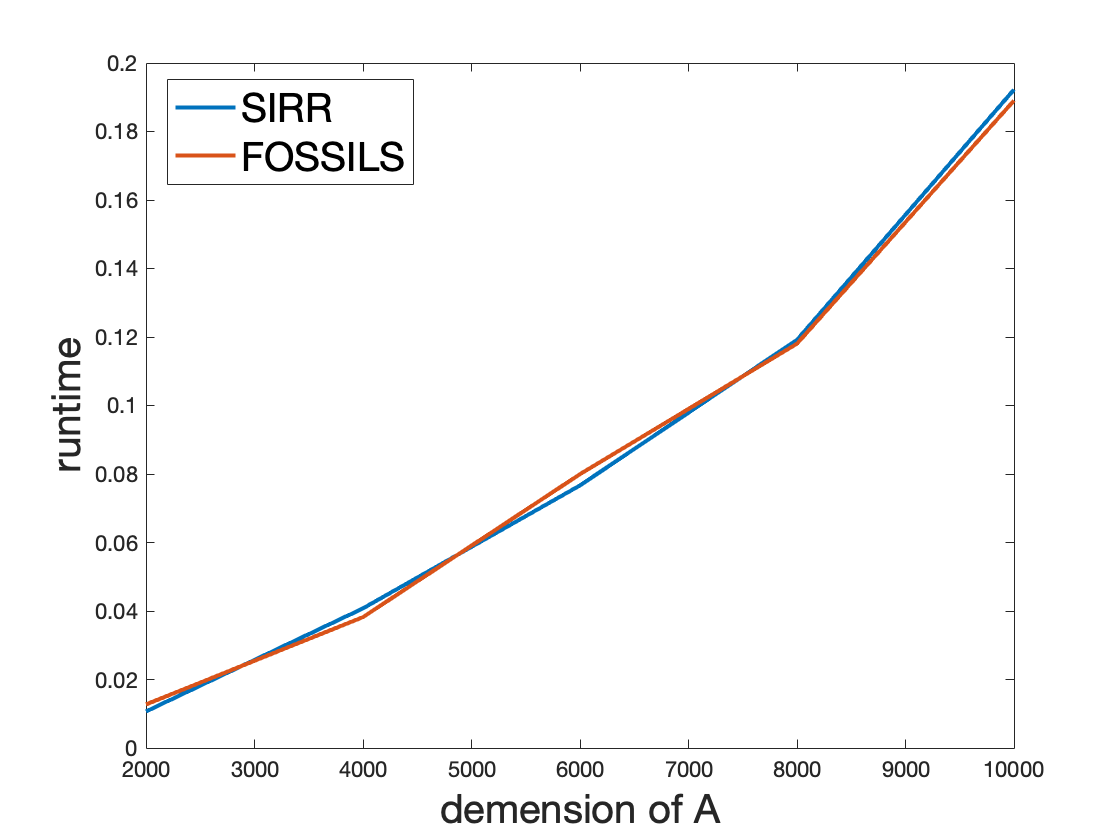}
    \includegraphics[width=2in]{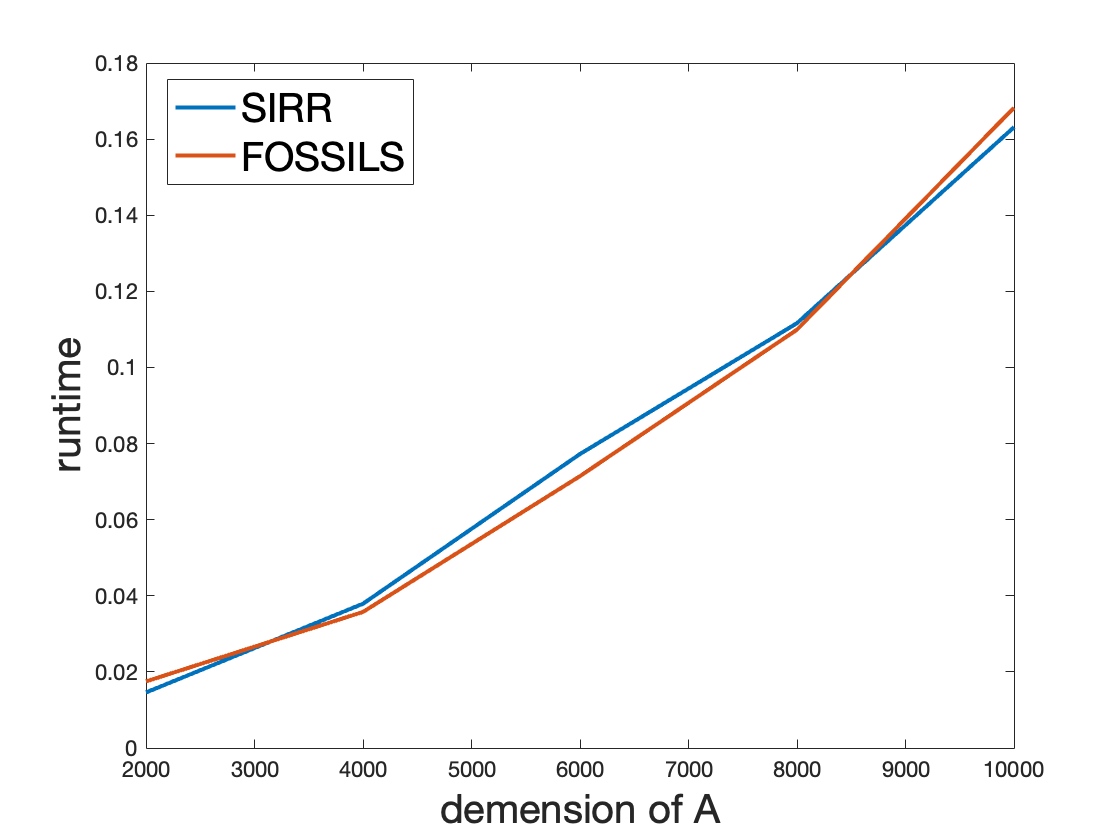}
    \caption{In 3 figures, $\kappa=10^4,10^8,10^{12}$ from left to right. SIRR has roughly the same computational cost as FOSSILS as randomized solver.}
    \label{fig:kerneltask} 
\end{figure*}

In the second experiment we give further insight into the difference between SIRR and FOSSILS.  First we test on kernel regression task to see the runtime of sketching solver and MATLAB solver (\textit{mldivide}) in different sizes of $n$. then we test the dependence of the stability of different solvers on sketching dimension by changing sketching dimension and counting the times that algorithm fails to converge in 100 runs. The left of Figure \ref{fig:kerneltask} shows that SIRR and FOSSILS need comparable time to reach the same accuracy, faster than MATLAB solver when $log(n)\geq 2.4$. The right figure illustrates the fail rate of SIRR and FOSSILS, which is the ratio of times failing to converge to a backward stable result in 100 runs. The fail rate of FOSSILS linearly decreases with the growth of sketching dimension, while SIRR achieves great stability when sketch dimension $d\geq 1.75 n$.

\begin{figure}[H]
    \centering
    \includegraphics[width=1.5in]{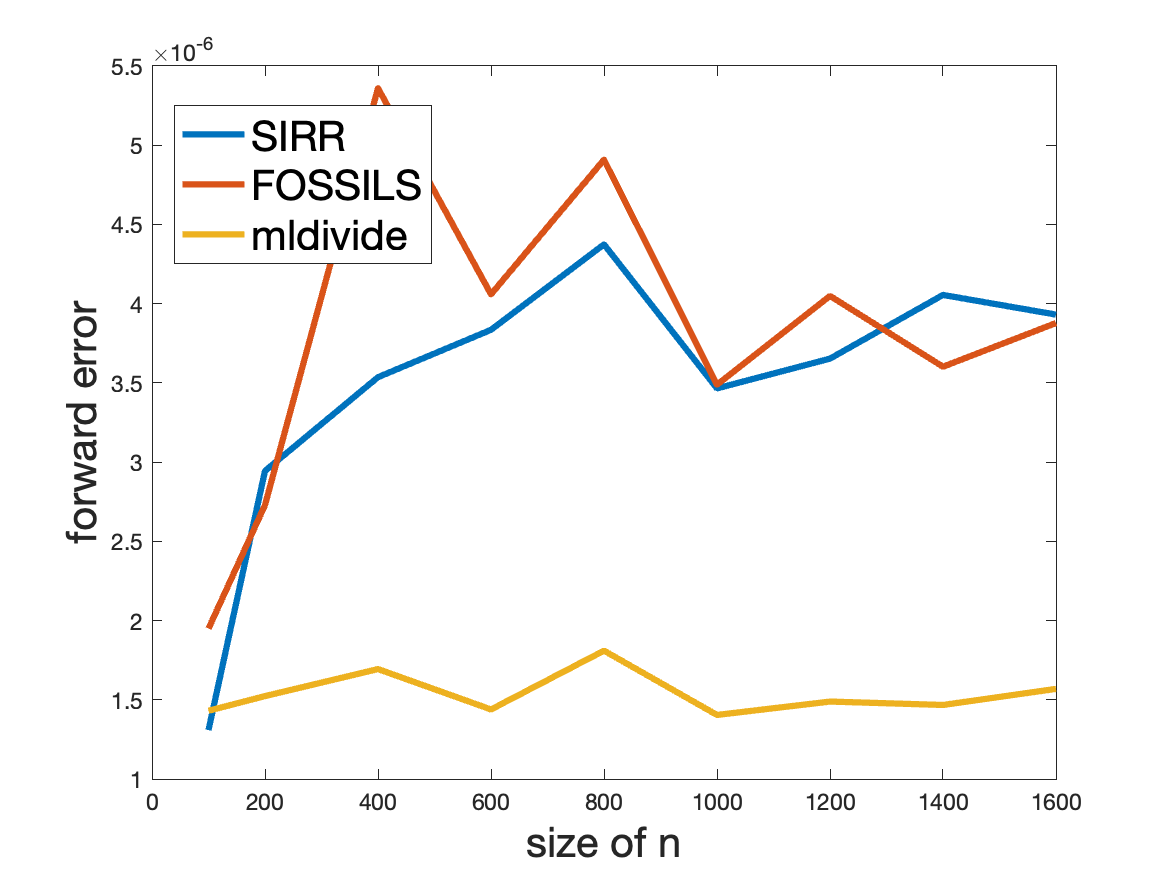}
    \includegraphics[width=1.5in]{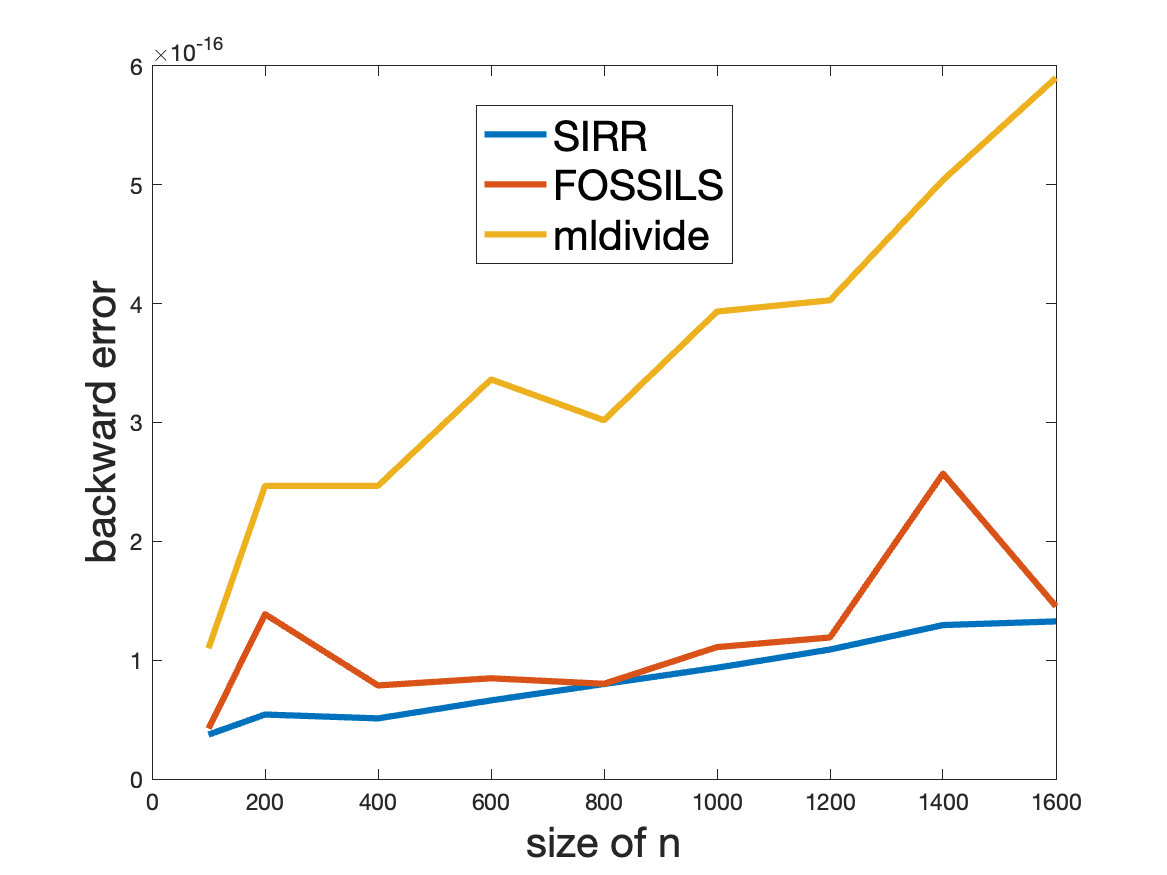}
   \caption{Forward error (left) and backward error (right) of different sizes of $n$.}
\label{fig:nscale}

\end{figure}

\paragraph{Error scale with $n$} In this experiment we show that for sketching solver and MATLAB direct solver, it is inevitable that the error is in scale with $n$. We fix $m=10000,\kappa=10^8,\|b-Ax^\ast\|=10^{-3}$ and change $n\in [100,1600]$ with sketching dimension $s=4*n$. Figure \ref{fig:nscale} shows the dependence of forward error and backward error on $n$ of different solvers. Note that three solvers actually have comparable forward error around $10^{-6}$ where MATLAB solver has slight edge. The dependence on $n$ is significant for backward error, where SIRR and FOSSILS appear to have a lower order of dependence. 
\begin{remark} Empirically, the growth of backward error as the matrix size increases is slower than the theoretical prediction in Remark \ref{nremark} as $n^2$. One possible reason is that the test matrix is random, and its randomness may not behave adversarially, leading to better performance. 
\end{remark}

\subsection{Comparison with FOSSILS}

We would like to highlight a concurrent work \cite{epperly2024fast}, which also developed a backward stable solver with a computational complexity of $O(mn+n^3)$. However, the FOSSILS solver proposed in their work follows a two-stage approach, where each stage involves an iterative process. In contrast, our algorithm is a single-stage solver that offers the flexibility to stop at any point during the computation, making it more adaptable for scenarios where early termination is necessary or beneficial. In this section we also compare our SIRR solver with the FOSSILS solver in both synthetic matrices (Figure \ref{fig:SIRRvsFOSSILS}) and realistic kernel regression datasets (Figure \ref{fig:nscale}). We demonstrate that the SIRR solver consistently achieves better backward stability than the FOSSILS solver across various difficulty levels, while requiring a similar amount of computing time. Notably, when the sketch dimension is small, SIRR is less prone to failure compared to FOSSILS.

\begin{figure*}
    \centering
    \vspace{-0.1in}  \includegraphics[width=2in]{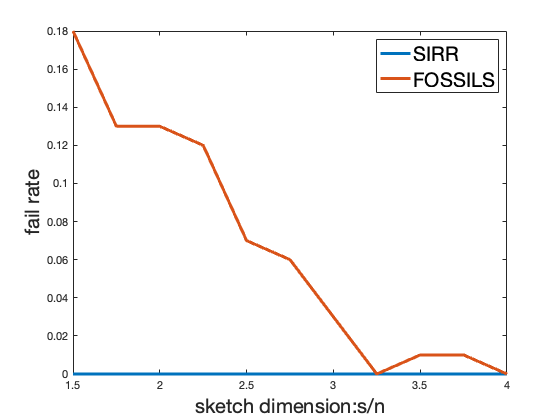}
    \includegraphics[width=2in]{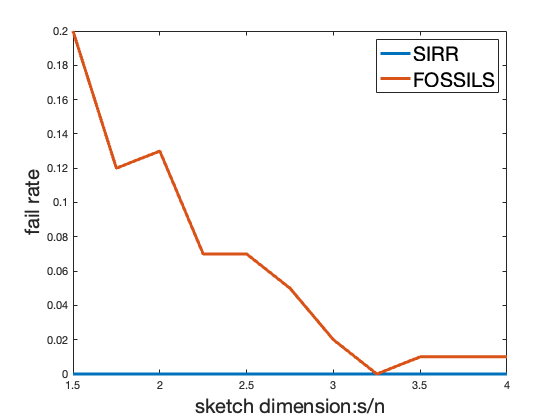}
    \includegraphics[width=2in]{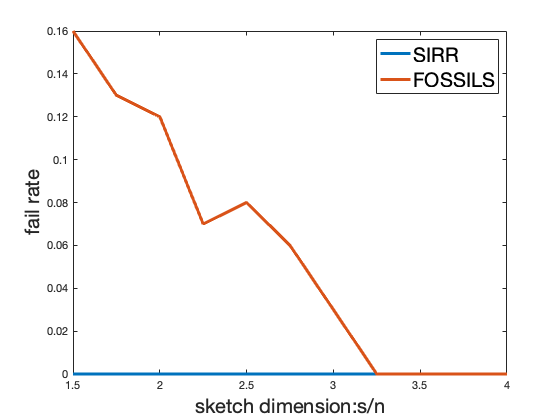}
    \includegraphics[width=2in]{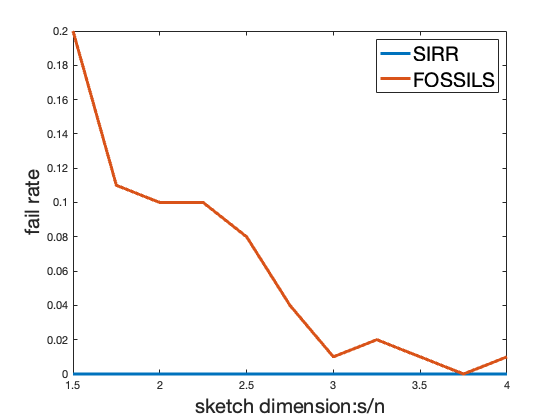}
    \includegraphics[width=2in]{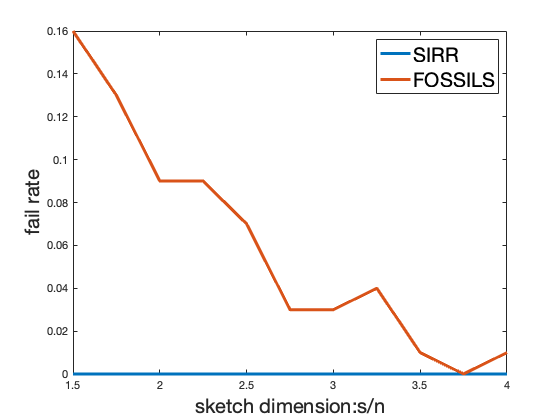}
    \includegraphics[width=2in]{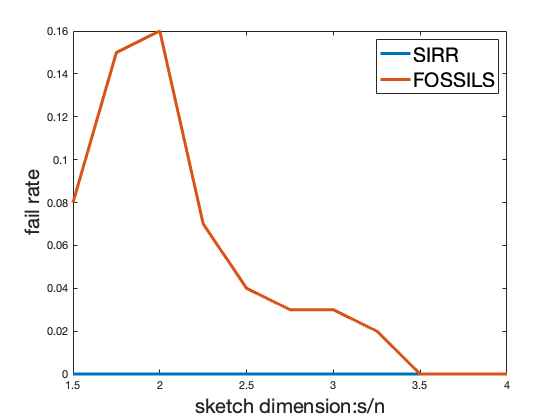}
    \caption{In first row $\|b-Ax^\ast\|=10^{-1}$ and in second row $\|b-Ax^\ast\|=10^{-3}$ with $\kappa = 10^4,10^8,10^{12}$ from left to right.}
    \label{dimension}
\end{figure*}

\paragraph{SIRR and FOSSILS with Different Embedding Quality}
 {We have proved that SIRR solver with 2-step Krylov-based meta-algorithm has a good convergence even in cases where sketching quality is bad and the distortion $\eta$ of sketching matrix is high. In this section, we give experiment results of the convergence performance of two solvers, SIRR and FOSSILS, in different embedding quality, which depends on the relative sketch dimensions $\frac{s}{n}$. In the result, the fail rate means the ratio of times that the solver fails to converge in 100 parallel experiments, where the criteria of failure in experiment is whether the relative residual error of result is under $10^{-5}$, given that the result in a good convergence has the relative residual error around $10^{-16}$. In different experiments, $m=2000$ and $n=100$, $\kappa\in\{10^4,10^8,10^{12}\}$, $\|b-Ax^\ast\|\in\{10^{-1},10^{-3}\}$. The results are presented in Figure \ref{dimension}.} 

\bibliographystyle{plain}
\bibliography{rho.bib}

\newpage
\onecolumn

\section*{Orgization of the Appendix}
The appendix is structured as follows. In section \ref{algo}, we present affiliate algorithms in this paper which are employed in practice. Then we give theoretical analysis of our method which guarantees the stability and convergence of our method in section \ref{proof}. 

\section{Details of Algorithms in Practice}
\subsection{random matrix}
 {In many applications, it is crucial to construct a subspace embedding without prior knowledge of the target subspace. Such embeddings are known as oblivious subspace embeddings. Typically, the singular values of specific random matrices are bounded with high probability, which makes them well-suited for subspace embedding. Various designs of random matrices exist that exhibit both strong computational and mathematical properties:
\begin{itemize}
    \item \textbf{Gaussian embedding}: $S=\mathbb{R}^{s\times m}$ with $i.i.d$ $ N(0,\frac{1}{s})$ entries. The normalization $\frac{1}{s}$ ensures that $S$ preserves the 2-norm in expectation, e.g. $E\|Sx\|_2^2=\|x\|_2^2$.
    \item \textbf{Subsampled randomized trigonometric transform (SRTT)}\cite{ailon2009fast}: $S=\sqrt{\frac{m}{s}}RDF\in \mathbb{R}^{s\times m}$, where $R\in \mathbb{R}^{s\times m}$ is an uniformly random set of $s$ rows drown from the identity matrix $I_m$, and $D\in \mathbb{R}^{m\times m}$ is a random diagonal matrix with $uniform(\pm 1)$ entries, and $F \in \mathbb{R}^{m\times m}$ is a $DCT_2$ matrix. SRTT requires less time in matrix/vector multiplication with cost $O(m\log(m))$ and has the same embedding property when $s\approx n\log(n)$.
    \item \textbf{Sparse random matrices}\cite{kaczmarz1937angenaherte}: $S=[s_1,s_2,\cdots,s_m]\in \mathbb{R}^{s\times m}$, where $s_i$ are sparse vectors, which means for each $i$, $s_i$ has exactly $\zeta$ nonzero entries, equally likely to be $\pm\sqrt{\frac{1}{\zeta}}$. The cost of matrix/vector multiplication is $O(\zeta m)$, and it's embedding has distortion $\eta$ when $s\approx \frac{n\log(n)}{\eta^2}$ and $\zeta\approx \frac{\log(n)}{\eta}$.
\end{itemize}
We use Sparse random matrices in our experiment which requires less operation and storage in computation.}
\subsection{Krylov-based meta-algorithm}
\label{algo}
 {A Krylov-based meta-algorithm is employed in our experiment, for it has a better convergence rate and is indifferent to the quality of embedding, making our solver more stable and faster even in worst cases. We present Krylov-based meta-algorithm in algorithm \ref{alg:krylov}}.
\begin{algorithm}
    \caption{Krylov-based meta-algorithm}\label{alg:krylov}
    \begin{tcolorbox}
    \hrule\hrule
    Krylov-based meta-algorithm 
    \hrule
    \SetKwInOut{Input}{Input}
    \SetKwInOut{Output}{Output}
    \SetKw{KwBy}{by}
    \SetKw{KwReturn}{Return}
    \SetKw{KwVia}{Via}
    \Input{1}
    \Output{2}
    \hrule\hrule
    $y_0= (R^\top R)^{-1}A^\top b$\;
    $Y=y_0$
    \Comment*[r]{Initialization via sketch-and-solve}
    \For{$i\gets1$ \KwTo $K$ \KwBy $1$}{
    $y_{i+1} = (R^\top R)^{-1}A^\top (b-Ay_i)$\;
    $Y = [Y,y+{i+1}]$
    }
    $a = (AY)^{\dag}b$
    
    \KwReturn $\text{ALG}^{\text{meta}}(b)=Ya$
    \end{tcolorbox}
    \end{algorithm}

\subsection{Sketched Iterative and Recursive Refinement}
 {Sketched Iterative and Recursive Refinement (SIRR) is provably fast and stable and is designed based on Sketched Iterative Refinement (SIR) and Sketched Recursive Refinement (SRR). We present SIRR in algorithm \ref{alg:SIRR}}.
\begin{algorithm}
\vspace{-0.1in}
\caption{\textbf{S}ketched \textbf{I}terative \textbf{R}ecursive \textbf{R}efinement.}\label{alg:SIRR}
\vspace{-0.1in}
\begin{tcolorbox}
\hrule\hrule
\textbf{SIRR}: Sketched Iterative Recursive Refinement
\hrule
\SetKwInOut{Input}{Input}
\SetKwInOut{Output}{Output}
\SetKw{KwBy}{by}
\SetKw{KwReturn}{Return}
\SetKw{KwVia}{Via}
\SetKw{KwInitial}{Initialize}
\Input{1}
\Output{2}
\hrule\hrule
    
    \If{$N=0$}{\KwInitial $\text{SIRR}_0(b) = (SA)^\dag Sb$ \;
    
    \Comment*[r]{Initialization via sketch-and-solve}}
    \For{$i\gets1$ \KwTo $N$ \KwBy $1$}{
    $\text{SIRR}_i(b):=\text{SIRR}_{i-1}(b)+{\color{orange}\text{SRR}_T} (A^\top (b-A \cdot \text{SIRR}_{i-1}(b)))$\Comment*[r]{Iterative Refinement via Sketched Recursive Refinement (SRR)}
    }
    \KwReturn $\text{SIRR}_N(b)$
    \vspace{-0.1in}
\end{tcolorbox}
\end{algorithm}

\section{PROOF OF MAIN RESULTS}
\label{proof}
In this section, we first establish some fundamental numerical results, which serve as the foundation for the subsequent numerical analysis. Then we examine the convergence of the iterative algorithm to show our method is theoretically fast. Finally, we give a rigorous numerical analysis of the algorithm to support that it is stable in both forward and backward sense.

\subsection{proof of lemma \ref{numerical error}}
 {
In this section, we prove some practical bounds for computed QR factorization $SA=\bar{Q}\hat{R}$.
The computation process of QR factorization can be decomposed as
\begin{align*}
    \widehat{SA} &= SA + E_1, \quad |E_1|\lesssim \gamma_n|S||A|,\\
    \widehat{SA}+E_2 &=\bar{Q}\hat{R} , \quad \|E_2\|_F\lesssim \gamma_{mn}\|\widehat{SA}\|_F,
\end{align*}
Thus we have
\begin{align*}
    \|\hat{R}\| = \|\bar{Q}\hat{R}\|=\|SA+E_1+E_2\|&\leq\frac{1}{1-\eta}\|A\|+2n\gamma_n\|A\|+2\sqrt{n}\gamma_{mn}\|A\|\lesssim \|A\|,\\
    \|\hat{R}^{-1}\| = \sigma_{min}(\bar{Q}\hat{R})=\sigma_{min}(SA+E_1+E_2)&\geq (1-\eta)\sigma_{min}(A)-(2\sqrt{n}\gamma_n\|A\|+2\sqrt{n}\gamma_{mn}\|A\|)\gtrsim \frac{\|A\|}{\kappa}.
\end{align*}
With similar analysis we have
\begin{align*}
    \|A\hat{R}^{-1}\|\leq \frac{1}{1-\eta}\|SA\hat{R}\|&\leq 2\|\bar{Q}-E_1\hat{R}^{-1}-E_2\hat{R}^{-1}\|\lesssim 1+u\kappa n^\frac{5}{2},\\
    \sigma_{min}(A\hat{R}^{-1})\geq (1-\eta)\sigma_{min}(SA\hat{R})&\geq\frac{1}{2}\sigma_{min}(\bar{Q}-E_1\hat{R}^{-1}-E_2\hat{R}^{-1})\gtrsim 1-u\kappa n^\frac{5}{2}.
\end{align*}
}
\subsection{proof of lemma \ref{becondition}}
 {In this section, we use straightforward computation to verify the relationship between $\alpha-\beta$ accuracy and backward error. In Karlson–Wald\'en estimate the key evaluating matrix can be expressed as}
\begin{align*}
&\left(A^{\top} A+\frac{\|b-A \widehat{x}\|^2}{1+\|\widehat{x}\|^2} I\right)^{-1 / 2}=\sum_{i=1}^n\left(\sigma_i^2+\frac{\|b-A \widehat{x}\|^2}{1+\|x\|^2}\right)^{-1 / 2} v_i v_i^{\top},
\end{align*}
 {
where $A = \sum_{i=1}^{n} \sigma_iu_iv_i^\top$ is SVD decomposition of matrix $A$. A further calculation shows that $\widehat{\mathrm{BE}}_1(\widehat{x})$ can be expressed as}
\begin{align*}
\widehat{\mathrm{BE}}_1(\widehat{x})^2&=\frac{1}{\sqrt{1+\|\widehat{x}\|^2}}\left\|\left(A^{\top} A+\frac{\|b-A \widehat{x}\|^2}{1+\|\widehat{x}\|^2} I\right)^{-1 / 2} A^{\top}(b-A \widehat{x})\right\|\\
&=\left\|\sum_{i=1}^n\left(\sigma_i^2+\frac{\|b-A \widehat{x}\|^2}{1+\|x\|^2}\right)^{-1 / 2} v_i v_i^{\top} A^{\top}A(x- \widehat{x})\right\|\\
&=\sum_{i=1}^n \frac{\sigma_i^4}{\left(1+\|\widehat{x}\|^2\right) \sigma_i^2+\|b-A \widehat{x}\|^2}\left(v_i^{\top}(x-\widehat{x})\right)^2.
\end{align*}
 {
Left multiplying \eqref{beform} by $v_i^{\top}$ yields }
\begin{align*}
    \left(v_i^{\top}(x-\widehat{x})\right)^2 \leq ((1+\|\widehat{x}\|)\sigma_i^{-1}+\|b-A \widehat{x}\|\sigma_i^{-2})^2.
\end{align*}

 {Combining two lines gives $\widehat{\mathrm{BE}}_1(\widehat{x})^2\leq n$.}

\subsection{the equivalence between $SIR$ and $SRR$}
\label{equivalence}
 {In this section, we show that sketched iterative refinement (SIR) and sketched recursive refinement (SIR) have the same form of results when they have the same linear meta-algorithm, which gives theoretical support to the statement that Recursive Refinement is just a reorganization of computation order
in the Iterative Refinement procedure.}

 {Set target linear system $Ax=b$ with $b\in range(A)$. Suppose that $ALG^{meta}(b) = Tb+q$ for some full rank matrix $T$, then we have }
\begin{align*}
    x_{i+1} = x_i + ALG^{meta}(b-Ax) = (I-TA)x + (Tb+q).
\end{align*}
 {
Note that the iteration is invariant when we initialize with true solution $x^\ast$ and $Ax^\ast=b$, thus $q=0$. A direct calculation shows that with a zero initial $x_0=0$, we have}
\begin{align*}
    \text{SIR}_N(b) = x_{N} = \sum_{i=0}^{N-1} (I-TA)^{i}Tb.
\end{align*}

 {For $x_i=\text{SRR}_i(b)$ and $x_0=0$, we claim that 
\begin{align*}
    \text{SRR}_N(b) = x_{N} = \sum_{i=0}^{2^N-1} (I-TA)^{i}Tb,
\end{align*}
and we will prove it by induction.}

 {It is easy to check that the statement holds for $N=1$ where
\begin{align*}
    \text{SRR}_0(b) &= ALG^{meta}(b) = Tb,\\
    \text{SRR}_1(b) & = ALG^{meta}(b)+ALG^{meta}(b-A\cdot ALG^{meta}(b))\\
    &=Tb+T(b-ATb)\\
    &=(I-TA)Tb+Tb.
\end{align*}}

 {To apply induction, suppose that the statement holds for $N$ and we compute $\text{SRR}_{N+1}(b)$ as
\begin{align*}
    \text{SRR}_{N+1}(b) & = \text{SRR}_{N}(b)+ALG^{meta}(b-A\cdot \text{SRR}_{N}(b))\\
    &=\sum_{i=0}^{2^N-1} (I-TA)^{i}Tb+\sum_{i=0}^{2^N-1} (I-TA)^{i}T(b-A\sum_{i=0}^{2^N-1} (I-TA)^{i}Tb).\\
\end{align*}
Note that $I-A\sum_{i=0}^{2^N-1} (I-TA)^{i}T = I-AT\sum_{i=0}^{2^N-1} (I-AT)^{i}=(I-AT)^{2^N}$, thus 
\begin{align*}
    \text{SRR}_{N+1}(b) &=\sum_{i=0}^{2^N-1} (I-TA)^{i}Tb+\sum_{i=0}^{2^N-1} (I-TA)^{i}T(b-A\sum_{i=0}^{2^N-1} (I-TA)^{i}Tb)\\
    &= \sum_{i=0}^{2^N-1} (I-TA)^{i}Tb + \sum_{i=0}^{2^N-1} (I-TA)^{i}T(I-AT)^{2^N}b\\
    &= \sum_{i=0}^{2^N-1} (I-TA)^{i}Tb + \sum_{i=2^N}^{2^{N+1}-1} (I-TA)^{i}Tb\\
    & =  \sum_{i=0}^{2^{N+1}-1} (I-TA)^{i}Tb.
\end{align*}
Thus the statement holds for all $N$.}

\subsection{proof of theorem \ref{convergence}}
\label{proofconvergence}
 {With a geometric series form of result given in the previous paper, one can easily examine the convergence of the iterative algorithm. Recall that for solving a well-defined linear system $Ax=b$, the solution of SIR and SRR has the form
\begin{align*}
     x_N = \sum_{i=0}^{N-1} (I-TA)^{i}Tb,
\end{align*}
and the true solution $x^\ast$ satisfies $x^\ast = A^{-1}b$. Thus the error $\|x_N-x^\ast\|$ satisfies
\begin{align*}
    \|x_N-x^\ast\| &= \|A^{-1}b-\sum_{i=0}^{N-1} (I-TA)^{i}Tb\|\\
    &= \|(I-\sum_{i=0}^{N-1} (I-TA)^{i}TA)A^{-1}b\|\\
    &= \|(I-TA)^{N}A^{-1}b\|\\
    &\leq \|(I-TA)\|^{N}\|A^{-1}b\|.
\end{align*}
It implies that SIR has a linear convergence with a convergence rate $\|I-TA\|$ and SRR has a quadratic convergence. The solver is convergent if and only if $\|I-TA\|< 1$.}

 {Then we compute the exact convergence rate for randomized solvers in solving $A^\top Ax = A^\top b$ instead of $Ax=b$, since in general cases $\|b-Ax^\ast\|\neq 0$. For iterative refinement with sketch-and-solve method, we have $T=(A^\top S^\top SA)^{-1}$ and thus $x_N = \sum_{i=0}^{N-1} (I-(A^\top S^\top SA)^{-1}A^\top A)^{i}(A^\top S^\top SA)^{-1}A^\top b$.}

 {Note that 
\begin{align*}
    \|(I-(A^\top S^\top SA)^{-1}A^\top A)^{i}(A^\top S^\top SA)^{-1}\| &= \|A^{\dag}(I-A(A^\top S^\top SA)^{-1}A^\top)^iA(A^\top S^\top SA)^{-1}\|\\
    &\leq \|A^{\dag}\|\|(I-A(A^\top S^\top SA)^{-1}A^\top)^i\|A(A^\top S^\top SA)^{-1}\|\\
    &\lesssim(\frac{\kappa}{\|A\|})^2(\frac{1}{(1-\eta)^2}-1)^i.
\end{align*}
The third inequality comes from the fact that $\|A(A^\top S^\top SA)^{-1}\|\asymp\|A^\dag\|=\frac{\kappa}{\|A\|}$ and $\|(I-A(A^\top S^\top SA)^{-1}A^\top)\|\leq \max\{\frac{1}{(1-\eta)^2}-1,1-\frac{1}{(1+\eta)^2}\}=\frac{1}{(1-\eta)^2}-1$. It implies
\begin{align*}
    \|x_N-x^\ast\| \leq constant \cdot(\frac{1}{(1-\eta)^2}-1)^i.
\end{align*}}

 {Note that to guarantee the convergence of SIR, the embedding distortion $\eta$ should be bounded in $(0,1)$. However, $\eta$ is usually bad in some difficult least-squares problems due to numerical error and small sketch dimensions. Fortunately, the Krylov subspace method is free from the restriction of $\eta$. We then verify the convergence of k-step Krylov-based iterative refinement. Note that 
\begin{align*}
    Ay_{i+1} &= A(I-(A^\top S^\top SA)^{-1}A^\top A)y_i+A (A^\top S^\top SA)^{-1} A^\top b\\
    &= (I-A(A^\top S^\top SA)^{-1}A^\top) Ay_i +  A(A^\top S^\top SA)^{-1} A^\top b\\
    y_0 &= x_i\\
    x_{i+1} &= \argmin_{x\in span\{y_1,y_2,\cdots,y_k\}} \|Ax-b\|.
\end{align*}
Denote $A(A^\top S^\top SA)^{-1}A^\top $ as $T$. Since $Ax_{i+1}\in Ax_i + \mathcal{K}_k (A(A^\top S^\top SA)^{-1}A^\top ,b-Ax_i)$, $A(x_{i+1}-x_i)$ can be expressed as 
$$A(x_{i+1}-x^\ast)=p_k(T)A(x_i-x^\ast),$$
where $p_k$ is a polynomial with order no more than $k$.}

 {Since $A(A^\top S^\top SA)^{-1}A^\top$ is normal matrix and can be decomposed as 
\begin{align*}
    A(A^\top S^\top SA)^{-1}A^\top = V\Lambda V^{\top},\quad V^\top V = I, \Lambda = diag(\lambda_1,\lambda_2,\cdots,\lambda_n),
\end{align*}
then
\begin{align*}
    A(x_{i+1}-x^\ast)=Vp_k(\Lambda)V^\top A(x_i-x^\ast),
\end{align*}
where $\|x_{i+1}-x^\ast\|_A$ is bounded by $|p_k(\lambda)|\|x_i-x^\ast\|$. We follow the practical but the worst case upper bound \cite{greenbaum1997iterative,liesen2004convergence} for min-max problem 
\begin{align*}
    \min_{p} \max_{i} p_k(\lambda_i)
\end{align*}
by choosing $p_k(\cdot)$ as the $k$-order Chebyshev polynomial. It leads to
\begin{align*}
    \|x_{i+1}-x^\ast\|_A &\leq (\frac{\sqrt{\kappa}-1}{\sqrt{\kappa+1}})^k\|x_i-x^\ast\|_A, \quad \kappa = \frac{\lambda_{max}(T)}{\lambda_{min}(T)},\\
\end{align*}
which further leads to
\begin{align*}
    \|x_{i+1}-x^\ast\|_A &\leq \min\{\eta^k,\frac{1}{\eta^k}\}\|x_i-x^\ast\|_A\\
\end{align*}
since $\kappa(T)=\kappa(A(A^\top S^\top SA)^{-1}A^\top)\leq \frac{(1+\eta)^2}{(1-\eta)^2}$. The result indicates that Krylov-based sketching method works even if the quality of subspace embedding is bad, requiring fewer sketching dimensions, which makes the algorithm faster.}

\subsection{proof of theorem \ref{forwardstable}}
 {In this section, we give a detailed analysis of the forward stability of SIRR, which also serves as a foundation for further discussion about backward stability.}

 {We first show the converged result of SRR can be decomposed into the form
$$\text{SRR}_N(r_A)\to x^\ast+u\sqrt{n}\|x^\ast\|e_1+u\kappa n^{\frac{3}{2}}\|Ax^\ast\|\hat{R}^{-1}e_2,$$
where $\|e_{1,2}\|\lesssim 1$.}

 {Consider the expression of $\text{SRR}_k(r_A)$ in real computation according to section \ref{equivalence}, we claim that 
\begin{align*}
   \hat{x}_k=\text{SRR}_{k}(r_A)& = (A^\top A)^{-1}r_A+a_k(\hat{R}^{-\top} r_A)e_k^1+b_k(\hat{R}^{-\top} r_A)R^{-1}e_k^2,\\
\end{align*}
where $a_k(\hat{R}^{-\top}r_A),b_k(\hat{R}^{-\top}r_A)$ are numerical errors, which are supposed to be functions of $\|\hat{R}^{-\top}r_A\|$, and $\|e_{i}^j\|\lesssim 1$. 
Then
\begin{align*}
    \text{SRR}_{k+1}(r_A)& = \hat{x}_k+ (A^\top A)^{-1}\hat{r}_k^A+a_k(R^{-\top}\hat{r}_k^A)e_k^5+b_k(R^{-\top}\hat{r}_k^A)R^{-1}e_k^6+\underbrace{\sqrt{n}\|\hat{x}_{k+1}\|}_{adding \ error}\\
    &=(A^\top A)^{-1}r_A+a_{k+1}(\hat{R}^{-\top} r_A)e_{k+1}^1+b_{k+1}(\hat{R}^{-\top} r_A)R^{-1}e_{k+1}^2,
\end{align*}
where
\begin{align*}
    \hat{r}_k^A&= r_A-A^\top A\hat{x}_k+un^\frac{3}{2}\|A\|\|\hat{x}_k\|A^\top e_k^3+un^\frac{3}{2}\|A\|\|A\hat{x}_k\|e_k^4.
\end{align*}
The first equation is a direct computation of $\text{SRR}_{k+1}$. The iteration of $a_k,b_k$ can thus be presented as
\begin{align}
\label{ab}
    b_{k+1}(\hat{R}^{-\top} r_A) &\lesssim b_k(\hat{R}^{-\top}\hat{r}_k^A)+un^\frac{3}{2}\|A\|\|\hat{x}_k\|+un^\frac{3}{2}\kappa \|A\hat{x}_k\|\nonumber\\
    a_{k+1}(\hat{R}^{-\top} r_A)&\lesssim a_k(\hat{R}^{-\top}\hat{r}_k^A) + u\sqrt{n}(\|\hat{x}-x^*\|+\|x^*\|)\nonumber
\end{align}
Note that
\begin{align*}
    R^{-\top }\hat{r}_k^A& = -a_k\|\hat{R}^{-\top} r_A\|R^{-\top }(A^\top A)e_k^1-b_k\|\hat{R}^{-\top} r_A\|R^{-\top }(A^\top A)R^{-1}e_k^2\\
    &+R^{-\top}(un^\frac{3}{2}\|A\|\|\hat{x}_k\|A^\top e_k^3+un^\frac{3}{2}\|A\|\|A\hat{x}_k\|e_k^4),\\
    \|R^{-\top }\hat{r}_k^A\|&\lesssim a_k\|A\|\|\hat{R}^{-\top} r_A\|+b_k\|\hat{R}^{-\top} r_A\|+u\kappa n^\frac{3}{2}\|\hat{R}^{-\top} r_A\|.
\end{align*}
where the last inequality comes from the fact that for $A^\top Ax_r^\ast=r$,
$$\frac{\|A\|}{\kappa}\|x^\ast_r\|\lesssim \|Ax^\ast_r\|\asymp\|\hat{R}^{-\top} r\|.
$$}

 {For $k=0$, the meta-algorithm is assumed to be $$\text{ALG}^{meta}(r_A) = (A^\top A)^{-1}r_A+c\|\hat{R}^{-\top}r_A\|\hat{R}^{-1}e$$
thus $a_0 = 0$ and $b_0 = c\|\hat{R}^{-1}r_A\|$. It's a natural idea to bound $a_k(\hat{R}^{-1}r_A)$ and $b_k(\hat{R}^{-1}r_A)$ by a linear function with respect to $\|\hat{R}^{-1}r_A\|$, since we can transform terms like $\|\hat{x}\|$, $\|A\hat{x}\|$ into $\|\hat{R}^{-1}r_A\|$ multiplied by some constant. First we convert terms $\|\hat{x}\|$ into $\|x^\ast\|+\|x^\ast-\hat{x}\|$ and convert terms $\|x^\ast\|$ and $\|Ax^\ast\|$ into $\|\hat{R}^{-\top}r_A\|$, and then compute $\hat{x}_k-x^\ast$ by leveraging the fact that $\hat{x}_k-x^\ast=a_k(\hat{R}^{-\top}r_A)e_k^1 + b_k(\hat{R}^{-\top}r_A) \hat{R}^{-1}e_k^2$. After assuming $a_k(\hat{R}^{-\top}r_A) \lesssim \alpha_k\|\hat{R}^{-\top}r_A\|$ and $b_k(\hat{R}^{-\top}r_A) \lesssim \beta_k\|\hat{R}^{-\top}r_A\|$, one gets the iteration of $\alpha_k, \beta_k$
\begin{align*}
    \beta_{k+1}\|\hat{R}^{-\top} r_A\| &\lesssim \beta_k(\alpha_k\|A\|\|\hat{R}^{-\top} r_A\|+\beta_k\|\hat{R}^{-\top} r_A\|+u\kappa n^\frac{3}{2}\|\hat{R}^{-\top} r_A\|)+un^\frac{3}{2}\|A\|\|\hat{x}_k\|+un^\frac{3}{2}\kappa \|A\hat{x}_k\|\\
    &\lesssim \alpha_k\beta_k\|A\|\|\hat{R}^{-\top} r_A\|+\beta_k^2\|\hat{R}^{-\top} r_A\|+u\kappa n^\frac{3}{2}\beta_k\|\hat{R}^{-\top} r_A\| + un^\frac{3}{2}\kappa \|\hat{R}^{-\top} r_A\|,\\
    a_{k+1}\|\hat{R}^{-\top} r_A\|&\lesssim \alpha_k(\alpha_k\|A\|\|\hat{R}^{-\top} r_A\|+\beta_k\|\hat{R}^{-\top} r_A\|+u\kappa n^\frac{3}{2}\|\hat{R}^{-\top} r_A\|) + u\sqrt{n}(\|\hat{x}-x^*\|+\|x^*\|)\\
    &\lesssim \alpha_k^2\|A\|\|\hat{R}^{-\top} r_A\|+\alpha_k\beta_k\|\hat{R}^{-\top} r_A\|+u\kappa n^\frac{3}{2}\alpha_k\|\hat{R}^{-\top} r_A\|,
\end{align*}
which leads to 
\begin{align*}
    \|A\|\alpha_{k+1}+\beta_{k+1} \lesssim ( \|A\|\alpha_{k}+\beta_{k})^2+un^{\frac{3}{2}}\kappa ( \|A\|\alpha_{k}+\beta_{k})+un^{\frac{3}{2}}\kappa.
\end{align*}
Since $\|A\|\alpha_{0}+\beta_{0}=\kappa n^{\frac{3}{2}}u<1$,  $\|A\|\alpha_{k}+\beta_{k}$ converges to $u\kappa n^{\frac{3}{2}}$ and thus $\alpha_k$ and $\beta_k$ converge to $u\kappa n^{\frac{3}{2}}$. Combined with \eqref{ab} one gets 
\begin{align*}
    \alpha_k&\lesssim (u\kappa n^{\frac{3}{2}})^2\|\hat{R}^{-\top}r_A\|+\sqrt{n}u\|x^\ast\|\lesssim \sqrt{n}u\|x^\ast\|,\\
    \beta_k&\lesssim un^{\frac{3}{2}}\|A\|\|x^\ast\|+u\kappa n^{\frac{3}{2}}\|Ax^\ast\|\lesssim u\kappa n^{\frac{3}{2}}\|Ax^\ast\|.
\end{align*}
Thus we have 
$$\text{SRR}(r_A)\to x^\ast+u\sqrt{n}\|x^\ast\|e_1+u\kappa n^{\frac{3}{2}}\|Ax^\ast\|\hat{R}^{-1}e_2.$$}

 {Now we iterate SRR with SIR to prove that SIRR is strongly forward stable. The real computation of SRR can be expressed as
\begin{align*}
    \hat{r}_i&=\underbrace{b-A\hat{x}_i}_{r_i}+u\sqrt{n}\|r_i\|e_{i,1}+un^{\frac{3}{2}}\|A\|\|\hat{x}_i\|e_{i,2},\\
    \hat{r}^A_i &= A^\top \hat{r}_i+un^{\frac{3}{2}}\|A\|\|\hat{r}_i\|e_{i,3},\\
    \hat{x}_{i+1} &= \hat{x}_i+(A^\top A)^{-1}(\hat{r}^A_i)+u\sqrt{n}\|x^\ast_r\|e_{i,4}+u\kappa n^{\frac{3}{2}}\|Ax^\ast_r\|\hat{R}^{-1}e_{i,5},\\
    &=(A^\top A)^{-1}A^\top b + (A^\top A)^{-1}(un^{\frac{3}{2}}\|A\|\|\hat{r}_i\|e_{i,3}+u\sqrt{n}\|r_i\|A^\top e_{i,1}+un^{\frac{3}{2}}\|A\|\|\hat{x}_i\|A^\top e_{i,2}),\\
    &+u\sqrt{n}\|x^\ast_r\|e_{i,4}+u\kappa n^{\frac{3}{2}}\|Ax^\ast_r\|\hat{R}^{-1}e_{i,5},\\
    x^\ast_r&=(A^\top A)^{-1}\hat{r}^A_i,\\
    &= (A^\top A)^{-1}(A^\top r_i+un^{\frac{3}{2}}\|A\|\|\hat{r}_i\|e_{i,3}+u\sqrt{n}\|r_i\|A^\top e_{i,1}+un^{\frac{3}{2}}\|A\|\|\hat{x}_i\|A^\top e_{i,2}),\\
\end{align*}
Denote $\|\hat{x}_i-x^\ast\|$ as $err_i$, $\|A\hat{x}_i-Ax^\ast\|$ as $err^r_i$. With decomposition $\|\hat{x}_i\|\leq\|x^\ast\|+err_i$, $\|A\hat{x}_i\|\leq\|Ax^\ast\|+err^r_i$, one gets
\begin{align*}
    \|x^\ast_r\|& \leq \frac{\kappa}{\|A\|}err_i^r+\frac{u\kappa^2 n^{\frac{3}{2}}}{\|A\|}\|r^\ast\|+u\kappa n^{\frac{3}{2}} \cdot err_i+u\kappa n^{\frac{3}{2}}\|x^\ast\|,\\
    \|Ax^\ast\|&\leq err_i^r + u\kappa n^{\frac{3}{2}} \|r^\ast\|+u n^{\frac{3}{2}}\|A\|err_i+un^{\frac{3}{2}}\|A\|\|x^\ast\|.
\end{align*}
We can then present the iteration of $err_i$ and $err^r_i$ as
\begin{align*}
    err_{i+1}&\leq \frac{u\kappa^2}{\|A\|}n^{\frac{3}{2}}(err^r_i+\|r^\ast\|)+u\kappa n^{\frac{3}{2}}(err_i+\|x^\ast\|)\\
    &+u\sqrt{n}(\frac{\kappa}{\|A\|}err_i^r+\frac{u\kappa^2}{\|A\|}\|r^\ast\|+u\kappa \cdot err_i+u\kappa\|x^\ast\|)\\
    &+\frac{u\kappa^2n^{\frac{3}{2}}}{\|A\|}(err_i^r + u\kappa\|r^\ast\|+u\|A\|err_i+u\|A\|\|x^\ast\|)\\
    &\leq u\kappa n^{\frac{3}{2}} \cdot err_i+\frac{u\kappa^2n^{\frac{3}{2}}}{\|A\|}err_i^r+\frac{u\kappa^2n^{\frac{3}{2}}}{\|A\|}\|r^\ast\|+u\kappa n^{\frac{3}{2}}\|x^\ast\|\\
    err^r_{i+1}&\leq u\kappa n^{\frac{3}{2}}(\|r^\ast\|+err_i^r)+un^{\frac{3}{2}}\|A\|(err_i+\|x^\ast\|)\\
    &+u\sqrt{n}\|A\|(\frac{\kappa}{\|A\|}err_i^r+\frac{u\kappa^2}{\|A\|}\|r^\ast\|+u\kappa \cdot err_i+u\kappa\|x^\ast\|)\\
    &+u\kappa n^{\frac{3}{2}}(err_i^r + u\kappa\|r^\ast\|+u\|A\|err_i+u\|A\|\|x^\ast\|)\\
    &=un^{\frac{3}{2}}\|A\|err_i+u\kappa n^{\frac{3}{2}}\cdot err_i^r+u\kappa n^{\frac{3}{2}}\|r^\ast\|+u n^{\frac{3}{2}}\|A\|\|x^\ast\|
\end{align*}
The iteration can be transformed into:
\begin{align*}
    \left(
        \begin{matrix}
            err_{i+1}\\
            err^r_{i+1}\\
            1
        \end{matrix}
    \right)&\lesssim n^{\frac{3}{2}}
    \left(
        \begin{matrix}
            u\kappa&\frac{u\kappa^2}{\|A\|}&\frac{u\kappa^2}{\|A\|}\|r^\ast\|+u\kappa\|x^\ast\|\\
            u\|A\|&u\kappa&u\|A\|\|x^\ast\|+u\kappa \|r^\ast\|\\
            0&0&1
        \end{matrix}
    \right)
    \left(
        \begin{matrix}
            err_{i}\\
            err^r_{i}\\
            1
        \end{matrix}
    \right)
\end{align*}
Since the transition matrix has the largest eigenvalue 1, the vector series $(err_{i+1},err^r_{i+1},1)^\top$ will converge to the eigenvector of 1, which leads to 
\begin{align*}
    \lim_{i\to \infty}
    \left(
        \begin{matrix}
            err_{i}\\
            err^r_{i}\\
            1
        \end{matrix}
    \right)&\lesssim n^{\frac{3}{2}}
    \left(
        \begin{matrix}
            u\kappa\|x^\ast\|+\frac{u\kappa^2}{\|A\|}\|r^\ast\|\\
            u\kappa\|r^\ast\|+u\|A\|\|x^\ast\|\\
            1
        \end{matrix}
    \right)\\
\end{align*}}

 {Thus the result of SIRR is forward stable. }
\subsection{proof of theorem \ref{strong solver}}
 { In this section, we propose the requirements for single step meta-solver to ensure that the SIR algorithm based on this meta-algorithm is backward stable. 
    Suppose that in $\text{i}^{th}$ iteration the solution $x_i$ has $a_i,b_i$-accuracy, which can be expressed as
    $$x_i=(A^\top A)^{-1}A^\top b+a_i\hat{R}^{-1}e_i^1+b_i(A^\top A)^{-1}e_i^2$$ for some unit random vector $e_i^1$ and $e_i^2$. We aim to get the iteration of $a_i, b_i$. Recall that $\tilde{\kappa}^{-1}=\max\{u\kappa, \frac{1}{\kappa n^{\frac{3}{2}}}\}$. 
    Following the computation of SIR we have
    \begin{align*}
        r_i &=b-Ax_i+f_i, \quad(\text{computed residual in each step})\\
        x_{i+1}&= x_i+(A^\top A)^{-1}(A^\top r_i)\\
        &+n^{\frac{3}{2}}(\tilde{\kappa}^{-1}\|Ax^\ast_r\|+u\kappa\tilde{\kappa}^{-1}\|r_i\|)\hat{R}^{-1}e_1\\
        &+n^{\frac{3}{2}}(u\|Ax^\ast_r\|+u\|r_i\|)(A^\top A)^{-1}e_2 \quad (\text{assumption of single step meta-solver})\\
        &=\underbrace{(A^\top A)^{-1}A^\top b}_{x^\ast} +a_{i+1}\hat{R}^{-1}e_{i+1}^1+b_{i+1}(A^\top A)^{-1}e_{i+1}^2,\\
    \end{align*}
    Here 
    \begin{align*}
        f_i &\lesssim u(n^{\frac{3}{2}}\|x_i\|+\sqrt{n}\|b-Ax_i\|)e_{f_i} ,\quad(\text{error in computed residual})\\
        x^\ast_r &=(A^\top A)^{-1}(A^\top r_i)\\
        &=-a_i(A^\top A)^{-1}A^\top (A \hat{R}^{-1})e_i^1-b_i(A^\top A)^{-1}e_i^2+(A^\top A)^{-1}A^\top f_i,\\
        b-Ax_i&=\underbrace{b-A(A^\top A)^{-1}A^\top b}_{r^\ast}-a_iA\hat{R}^{-1}e_i^1-b_iA(A^\top A)^{-1}e_i^2,\\
    \end{align*}
    We omit the adding error when computing $x_{i+1}=x_i+d_i$ where $d_i$ is the refinement term, since it is already machine-precision. Then we have
    \begin{align*}
        a_{i+1}\hat{R}^{-1}e_{i+1}^1&=n^{\frac{3}{2}}(\tilde{\kappa}^{-1}\|Ax^\ast_r\|+u\kappa\tilde{\kappa}^{-1}\|r_i\|)\hat{R}^{-1}e_1\\
        &+\underbrace{u(n^{\frac{3}{2}}\|x_i\|+\sqrt{n}\|b-Ax_i\|)\hat{R}^{-1}(\hat{R}(A^\top A)^{-1}A^\top e_{f_i})}_{(A^\top A)^{-1}A^\top f_i},\\
        b_{i+1}(A^\top A)^{-1}e_{i+1}^2&=n^{\frac{3}{2}}(u\|Ax^\ast_r\|+u\|r_i\|)(A^\top A)^{-1}e_2,
    \end{align*}
    which yields
    \begin{align*}
        a_{i+1}& \lesssim n^{\frac{3}{2}}(\tilde{\kappa}^{-1}\|Ax^\ast_r\|+u\kappa\tilde{\kappa}^{-1}\|r_i\|)+\|f_i\|\\
        &\lesssim n^{\frac{3}{2}}\tilde{\kappa}^{-1} (a_i+\kappa b_i+\|f_i\|)\\
        &+n^{\frac{3}{2}}u\kappa\tilde{\kappa}^{-1} (\|r^\ast\|+a_i+\kappa b_i+\|f_i\|)\\
        &+\underbrace{u(n^{\frac{3}{2}}\|x_i\|+\sqrt{n}(\|r^\ast\|+a_i+\kappa b_i))}_{\|f_i\|}\\
        &\asymp n^{\frac{3}{2}}\tilde{\kappa}^{-1} (a_i +\kappa b_i) + un^{\frac{3}{2}}\|x^\ast\| + u\kappa n^{\frac{3}{2}}\tilde{\kappa}^{-1}\|r^\ast\|,\\
        b_{i+1}&\lesssim un^{\frac{3}{2}}\|Ax^\ast_r\|+un^{\frac{3}{2}}\|r_i\|\\
        &\lesssim un^{\frac{3}{2}}(a_i+\kappa b_i+\|f_i\|)\\
        &+un^{\frac{3}{2}}(\|r^\ast\|+a_i+\kappa b_i+\|f_i\|)\\
        &\asymp un^{\frac{3}{2}}(a_i+\kappa b_i)+u^2n^3\|x^\ast\|+un^{\frac{3}{2}}\|r^\ast\|.
    \end{align*}}
    
     {The iteration of $a_i,b_i$ can be written in the form
    \begin{align*}
        \left(
            \begin{matrix}
                a_{i+1}\\
                b_{i+1}\\
                1
            \end{matrix}
        \right)&\lesssim
        \left(
            \begin{matrix}
                n^{\frac{3}{2}}\tilde{\kappa}^{-1} &n^{\frac{3}{2}}\kappa\tilde{\kappa}^{-1}&un^{\frac{3}{2}}\kappa\tilde{\kappa}^{-1} \|r^\ast\|+un^{\frac{3}{2}}\|x^\ast\|\\
                un^{\frac{3}{2}}&un^{\frac{3}{2}}\kappa&un^{\frac{3}{2}}\|r^\ast\|+u^2n^3\|x^\ast\|\\
                0&0&1
            \end{matrix}
        \right)
        \left(
            \begin{matrix}
                a_{i}\\
                b_{i}\\
                1
            \end{matrix}
        \right)\\
    \end{align*}
    Since the transition matrix has the largest eigenvalue 1, the vector series $(a_i,b_i,1)^\top$ will converge to the eigenvector of 1, which leads to 
    \begin{align*}
        \lim_{i\to \infty}
        \left(
            \begin{matrix}
                a_{i}\\
                b_{i}\\
                1
            \end{matrix}
        \right)&\lesssim
        \left(
            \begin{matrix}
                un^3\kappa\tilde{\kappa}^{-1}\|r^\ast\|+un^{\frac{3}{2}}\|x^\ast\|\\
                un^{\frac{3}{2}}\|r^\ast\|+u^2n^3\|x^\ast\|\\
                1
            \end{matrix}
        \right),\\
    \end{align*}
    where we use the fact that $\tilde{\kappa}^{-1}n^{\frac{3}{2}}<1$.}

     {With $u\kappa\tilde{\kappa}^{-1}=u\kappa\max\{u\kappa,\kappa^{-1}\}<u^2\kappa^2+u$, the result $\lim_{i\to\infty}x_i$ is backward stable.}
    \subsection{proof of theorem \ref{iterativeIDS}}
     { In this section, we show that SRR generates a good single step meta-solver for SIR, in other words, SIRR is backward stable. In SRR, the meta-algorithm $\text{SRR}_0(\cdot)$ solves a full rank linear system $(A^\top A)x=r_A$, and from iteration process we find the error of solution only depends on either $\|R^{-\top }r_A\| $, $\|x^\ast\|$ or $\|Ax^\ast\|$, where $x^\ast = (A^\top A)^{-1}r_A$ and $\|R^{\-top} r_A\|\asymp \|Ax^\ast\|$. Thus we can assume that 
    $$\text{SRR}_i(r_A) = (A^\top A)^{-1}r_A +  (a_i^1 \|x^\ast\| + a_i^2 \|Ax^\ast\|)R^{-1}e_1^i + (b_i^1 \|x^\ast\| + b_i^2 \|Ax^\ast\|)(A^\top A)^{-1}e_2^i
    $$
    , one can get the iteration of $a_i^j, b_i^j$ w.r.t $i$.}

     {The iteration of $\text{SRR}_i$ then can be written as
    \begin{align*}
        \text{SRR}_{i+1}(r_A) &= \text{SRR}_i(r_A) + \text{SRR}_i(\underbrace{r_A-A^\top A x_i+f_i}_{\hat{r}_i})\\
         &=x_i + (A^\top A)^{-1}(r_A-A^\top A x_i)\\
        &+(A^\top A)^{-1}(\underbrace{(u\sqrt{n}\|r_A-A^\top A x_i\|+un^{\frac{3}{2}}\|Ax_i\|)e_{f_1}+un^{\frac{3}{2}}\|x_i\|A^\top e_{f_2})}_{f_i}\\
        &+(a_i^1 \|x^\ast_r\| + a_i^2 \|Ax^\ast_r\|)\hat{R}^{-1}e_1^i + (b_i^1 \|x^\ast_r\| + b_i^2 \|Ax^\ast_r\|)(A^\top A)^{-1}e_2^i\\
        &=(A^\top A)^{-1}r_A +(a_{i+1}^1 \|x^\ast_r\| + a_{i+1}^2 \|Ax^\ast_r\|)\hat{R}^{-1}e_1^{i+1} \\
        &+ (b_{i+1}^1 \|x^\ast_r\| + b_{i+1}^2 \|Ax^\ast_r\|)(A^\top A)^{-1}e_2^{i+1},\\
    \end{align*}
    where
    \begin{align*}
        x^\ast_r & = (A^\top A)^{-1}\hat{r}_i\\
        & = (A^\top A)^{-1}(r_A - A^\top Ax_i+f_i)\\
        & = -(a_i^1 \|x^\ast\| + a_i^2 \|Ax^\ast\|)\hat{R}^{-1}e_1^i - (b_i^1 \|x^\ast\| + b_i^2 \|Ax^\ast\|)(A^\top A)^{-1}e_2^i + (A^\top A)^{-1}f_i,\\
        f_i& = u\sqrt{n}\|r_A-A^\top A x_i\|+un^{\frac{3}{2}}\|Ax_i\|)e_{f_1}+un^{\frac{3}{2}}\|x_i\|A^\top e_{f_2}.\\
    \end{align*}
    Denote $1+u\kappa n^{\frac{3}{2}}$ as $\hat{1}$ for convenience, then the expansion of $\|x^\ast_r\|$ and $\|Ax^\ast_r\|$ yields
    \begin{align*}
        \|x^\ast_r\|&\lesssim (\kappa a_i^1 +\kappa^2 b_i^1)\|x^\ast\| + (\kappa a_i^2+ \kappa^2 b_i^2)\|Ax^\ast\| \\
        &+\underbrace{u\kappa^2\sqrt{n}((a_i^1+b_i^1)\|x^\ast\|+(a_i^2+b_i^2)\|Ax^\ast\|)}_{f_i \ term1}\\
        &+\underbrace{u\kappa^2n^{\frac{3}{2}}(\|Ax^\ast\|+(a_i^1+\kappa b_i^1)\|x^\ast\|+(a_i^2+\kappa b_i^2)\|Ax^\ast\|)}_{f_i \ term2}\\
        &+\underbrace{u\kappa n^{\frac{3}{2}}(\|x^\ast\|+\kappa(a_i^1+\kappa b_i^1)\|x^\ast\|+\kappa(a_i^2+\kappa b_i^2)\|Ax^\ast\|)}_{f_i \ term3}\\
        &\asymp (\hat{1}\kappa a_i^1 + \hat{1}\kappa^2 b_i^1 + u\kappa n^{\frac{3}{2}})\|x^\ast\|+(\hat{1}\kappa a_i^2+\hat{1}\kappa^2 b_i^2+u\kappa^2n^{\frac{3}{2}})\|Ax^\ast\|,\\
    \end{align*}       
        \begin{align*}
        \|Ax^\ast_r\|&\lesssim ( a_i^1 +\kappa b_i^1)\|x^\ast\| + (a_i^2+ \kappa b_i^2)\|Ax^\ast\| \\
        &+\underbrace{u\kappa\sqrt{n}((a_i^1+b_i^1)\|x^\ast\|+(a_i^2+b_i^2)\|Ax^\ast\|)}_{f_i \ term1}\\
        &+\underbrace{u\kappa n^{\frac{3}{2}}(\|Ax^\ast\|+(a_i^1+\kappa b_i^1)\|x^\ast\|+(a_i^2+\kappa b_i^2)\|Ax^\ast\|)}_{f_i \ term2}\\
        &+\underbrace{u n^{\frac{3}{2}}(\|x^\ast\|+\kappa(a_i^1+\kappa b_i^1)\|x^\ast\|+\kappa(a_i^2+\kappa b_i^2)\|Ax^\ast\|)}_{f_i \ term3}\\
        &\asymp (\hat{1}a_i^1 + \hat{1}\kappa b_i^1 + un^{\frac{3}{2}})\|x^\ast\|+(\hat{1} a_i^2+\hat{1}\kappa b_i^2+u\kappa n^{\frac{3}{2}})\|Ax^\ast\|,\\
    \end{align*}}
    
     {With assumption $u\kappa n^{\frac{3}{2}}<1$, $\hat{1}\lesssim 1$, the iteration of $a_i^j,b_i^j$ has the form
    \begin{align*}
        \left(
            \begin{matrix}
                a_{i+1}^1\\
                a_{i+1}^2\\
            \end{matrix}
        \right)&\lesssim 
        \left(
            \begin{matrix}
                \kappa (a_i^1 + \kappa b_i^1+un^{\frac{3}{2}}) & (a_i^1 + \kappa b_i^1+un^{\frac{3}{2}})\\
                \kappa(a_i^2+\kappa b_i^2+u\kappa n^{\frac{3}{2}}) & (a_i^2+\kappa b_i^2+u\kappa n^{\frac{3}{2}})
            \end{matrix}
        \right)
        \left(
            \begin{matrix}
                a_{i}^1\\
                a_{i}^2\\
            \end{matrix}
        \right)
        +n^{\frac{3}{2}}
        \left(
            \begin{matrix}
                u+u\kappa a_i^1+u\kappa^2 b_i^1\\
                u\kappa a_i^2+u\kappa^2 b_i^2\\
            \end{matrix}
        \right),\\
        \left(
            \begin{matrix}
                b_{i+1}^1\\
                b_{i+1}^2\\
            \end{matrix}
        \right)&\lesssim 
        \left(
            \begin{matrix}
                \kappa (a_i^1 + \kappa b_i^1+un^{\frac{3}{2}}) & (a_i^1 + \kappa b_i^1+un^{\frac{3}{2}})\\
                \kappa(a_i^2+\kappa b_i^2+u\kappa n^{\frac{3}{2}}) & (a_i^2+\kappa b_i^2+u\kappa n^{\frac{3}{2}})
            \end{matrix}
        \right)
        \left(
            \begin{matrix}
                b_{i}^1\\
                b_{i}^2\\
            \end{matrix}
        \right)
        +n^{\frac{3}{2}}
        \left(
            \begin{matrix}
                u(a_i^1+\kappa b_i^1)\\
                u(a_i^2+\kappa b_i^2+1)\\
            \end{matrix}
        \right).\\
    \end{align*}}
    
     {Let $c_i^j = a_i^j+\kappa b_i^j$, then}
    \begin{align*}
        c_{i+1}^1 &\lesssim \kappa (c_i^1)^2+c_i^1c_i^2+ u\kappa n^{\frac{3}{2}} c_i^1 + un^{\frac{3}{2}}c_i^2+un^{\frac{3}{2}},\\
        c_{i+1}^2 &\lesssim \kappa c_i^1 c_i^2 + (c_i^2)^2+ u\kappa^2 n^{\frac{3}{2}}c_i^1 + u\kappa n^{\frac{3}{2}}c_i^2 + u\kappa n^{\frac{3}{2}},\\
        \kappa c_{i+1}^1+c_{i+1}^2&\lesssim(\kappa c_i^1+c_i^2)^2+u\kappa n^{\frac{3}{2}}(\kappa c_i^1+c_i^2)+u\kappa n^{\frac{3}{2}}.\\
    \end{align*}
    
     {For general cases, consider the initialization of 
    \begin{align*}
        \kappa a_0^1+a_0^2\asymp c, \kappa(\kappa b_0^1+b_0^2)\asymp c, c_0\asymp\kappa c_0^1+c_0^2\asymp \kappa(a_0^1+\kappa b_0^1) + a_0^2+\kappa b_0^2\asymp c.
    \end{align*}
    Then $\kappa c_0^1+c_0^2=c\gtrsim u\kappa n^{\frac{3}{2}}$, thus
    \begin{align*}
        \kappa c_i^1+c_i^2 &\lesssim  max\{(\underbrace{\kappa c_0^1+c_0^2}_{c})^{2^i},u\kappa n^{\frac{3}{2}}\}.\\
    \end{align*}
    Similar result can be derived for $c_i^1,c_i^2$ where
    \begin{align*}
        c_{i+1}^1 &\lesssim (\kappa c_i^1+c_i^2)c_i^1+ un^{\frac{3}{2}}(\kappa c_i^1+c_i^2)+un^{\frac{3}{2}}\lesssim max\{ (\kappa c_i^1+c_i^2)c_i^1,u n^{\frac{3}{2}}\},\\
        c_{i+1}^2 &\lesssim (\kappa c_i^1+c_i^2)c_i^2 + u\kappa n^{\frac{3}{2}}(\kappa c_i^1+c_i^2) + u\kappa n^{\frac{3}{2}}\lesssim max\{(\kappa c_i^1+c_i^2) c_i^2,u\kappa n^{\frac{3}{2}}\}\,\
    \end{align*}
    which leads to
    \begin{align*}
        c_{i}^1&\lesssim max\{c_0^1(c_0)^{2^i},u n^{\frac{3}{2}}\},\\
        c_{i}^2&\lesssim max\{c_0^2(c_0)^{2^i},u\kappa n^{\frac{3}{2}}\}.\\
    \end{align*}}
    
     {For $a_i^j$, we can first calculate the iteration of $\kappa a_i^1+a_i^2$by 
    \begin{align*}
        \kappa a_{i+1}^1+a_{i+1}^2&\lesssim (\underbrace{\kappa c_i^1+c_i^2}_{\lesssim max\{(\kappa c_0^1+c_0^2)(c)^{2^i},u\kappa n^{\frac{3}{2}}\}}) (\kappa a_i^1+a_i^2) + u\kappa n^{\frac{3}{2}}(\kappa a_i^1+a_i^2) \\
        &+ u\kappa n^{\frac{3}{2}} +u\kappa n^{\frac{3}{2}}(\kappa c_i^1+c_i^2) \quad(\text{from transition matrix})\\
        &\lesssim (\kappa c_i^1+c_i^2)(\kappa a_i^1+a_i^2)+u\kappa n^{\frac{3}{2}}\quad(\sup_{i} (\kappa c_i^1+c_i^2)>u\kappa n^{\frac{3}{2}},\sup_{i}\kappa a_i^1+a_i^2>u\kappa n^{\frac{3}{2}})\\
        &\lesssim max\{(c)^{2^{i+1}}(\kappa a_0^1+a_0^2),u\kappa n^{\frac{3}{2}}\},
    \end{align*}
    thus
    \begin{align*}
        a_{i+1}^1&\lesssim (\kappa a_i^1+a_i^2)(c_i^1+u n^{\frac{3}{2}})+u n^{\frac{3}{2}}+u\kappa n^{\frac{3}{2}}c_i^1\\
        &\lesssim max\{(c)^{2^{i+1}}(\kappa a_0^1+a_0^2)c_0^1,un^{\frac{3}{2}}\}.
    \end{align*}
    Similarly
    \begin{align*}
        a_i^2&\lesssim max\{(c)^{2^{i}}(\kappa a_0^1+a_0^2)c_0^2,u^2\kappa^2n^3\},\\
        b_i^1&\lesssim max\{(c)^{2^{i}}(\kappa b_0^1+b_0^2)c_0^1,u^2n^3\},\\
        b_i^2&\lesssim max\{(c)^{2^{i}}(\kappa b_0^1+b_0^2)c_0^2,un^{\frac{3}{2}}\},\\
    \end{align*}
    which leads to the bound of $\text{SRR}_\infty(r_A)$:
    \begin{align*}
        \text{SRR}_\infty(r_A) &= \underbrace{(A^\top A)^{-1}r_A}_{x^\ast} + \hat{a}\hat{R}^{-1}e_1+\hat{b}(A^\top A)^{-1}e_2
    \end{align*}
    where 
    $$\hat{a}= un^{\frac{3}{2}}\|x^\ast\|+u^2\kappa^2 n^3\|Ax^\ast\|,\hat{b}= u^2n^3\|x^\ast\|+un^{\frac{3}{2}}\|Ax^\ast\|, \quad \|e_{1,2}\|\lesssim 1.
    $$}
    
     {In iteration algorithm, we need to compute $A^\top b$ as $r_A$ with $error(r_A)=un^{\frac{3}{2}}\|b\|e$, so $r_A=A^\top b+u n^{\frac{3}{2}}\|b\|e$. Then $\text{SRR}_\infty$ becomes
    \begin{align*}
        \text{SRR}_\infty(b) &= (A^\top A)^{-1}A^\top b + a\hat{R}^{-1}e_1+b(A^\top A)^{-1}e_2,
    \end{align*}
    where 
    \begin{align*}
        a&\asymp un^{\frac{3}{2}}\|\hat{x}^\ast\|+u^2\kappa^2n^3\|A\hat{x}^\ast\|\\
        &=un^{\frac{3}{2}}\|x^\ast+un^{\frac{3}{2}}\|b\|(A^\top A)^{-1}e\|+u^2\kappa^2n^3\|Ax^\ast+un^{\frac{3}{2}}\|b\|A(A^\top A)^{-1}e\|\\
        &=un^{\frac{3}{2}}\|x^\ast\|+u^2\kappa^2n^3\|Ax^\ast\|+u^2\kappa^2n^3\|b\|,\\
        b&\asymp u^2n^3\|\hat{x}^\ast\|+un^{\frac{3}{2}}\|A\hat{x}^\ast\|+un^{\frac{3}{2}}\|b\|\\
        &=u^2n^3\|x^\ast\|+un^{\frac{3}{2}}\|Ax^\ast\|+un^{\frac{3}{2}}\|b\|.
    \end{align*}}
    
     {However, for practical use, we can stop the iteration as soon as the algorithm achieves the accuracy needed in theorem \ref{strong solver}.For $\text{SRR}_N$ with Two-step Krylov-based meta-algorithm, $a_0^1=b_0^1=0,a_0^2\asymp \kappa b_0^2\asymp c$, the iteration only refines $a_i^2$ and $b_i^2$. With $(\kappa a_0^1+a_0^2)\asymp 1, (\kappa b_0^1+b_0^2)\kappa\asymp 1$, the steps we actually need in practice is 
    $N=\max\{log_2(\frac{log(\tilde{\kappa }^{-1}n^{\frac{3}{2}})}{log(c)}),log_2(\frac{log(u\kappa n^{\frac{3}{2}})}{log(c)})\} = log_2(\frac{log(\tilde{\kappa }^{-1}n^{\frac{3}{2}})}{log(c)})$, with
    \begin{align*}
        \text{SRR}_N(b) &= (A^\top A)^{-1}A^\top b + n^{\frac{3}{2}}(u\|x^\ast\|+\tilde{\kappa}^{-1}\|Ax^\ast\|)\hat{R}^{-1}e_1+(u^2n^3\|x^\ast\|+un^{\frac{3}{2}}\|Ax^\ast\|)(A^\top A)^{-1}e_2.
    \end{align*}
    A similar discussion further leads to the result of considering the first multiplication $r_A=A^\top b$
    \begin{align*}
        \text{SRR}_N(b) &= (A^\top A)^{-1}A^\top b \\
        &+n^{\frac{3}{2}}(\tilde{\kappa}^{-1}\|Ax^\ast\|+u\|x^\ast\|+u\kappa\tilde{\kappa}^{-1}\|b\|)\hat{R}^{-1}e_1\\
        &+n^{\frac{3}{2}}(u\|Ax^\ast\|+u^2n^{\frac{3}{2}}\|x^\ast\|+u\|b\|)(A^\top A)^{-1}e_2\\
    \end{align*}}

\subsection{Proof of Lemma \ref{metasolver}}
 {
In this section, we verify that the Krylov-based meta solver satisfies the condition of theorem \ref{iterativeIDS}, which finally proves that a k-step Krylov-based SIRR solver is backward stable.} {In Krylov subspace method, with $y_0,y_1,\cdots, y_k$ given by using iterative sketching, we solve the least squares problem in the space spanned by $\{y_i\}_{i=1}^k$
$$\argmin_{x\in span\{y_0,y_1,\cdots, y_k\}}\|r_A-(A^\top A)x\|.$$
Let $Y:=[y_0,y_1,\cdots,y_k]$, then the solution is $x=Y(A^\top AY)^{-1}r_A$. }

 {Consider the numerical process of computing $Y(A^\top AY)^{-1}r_A$, which is
\begin{align*}
    \widehat{AY} &= AY+E_1,\quad\|E_1\|\lesssim un^{\frac{3}{2}}\|A\|\|Y\|,\\
    (A^\top \widehat{AY}+E_2)\hat{a}&=  r_A+ h_1,\quad \|E_2\|\leq un^{\frac{3}{2}}\sqrt{k}\|A\|\|AY\|+uk^{\frac{5}{2}}\|A^\top \widehat{AY}\|,\|h_1\|\leq uk^2\|r_A\|,\\
    \hat{x}&=Y\hat{a}+h_2,\quad \|h_2\|\lesssim uk^\frac{3}{2}\|Y\|\|\hat{a}\|.
\end{align*}}

 {Note that $y_i$ can be expressed as $y_i=(A^\top A)^{-1}r_A+c\|\hat{R}^{-\top} r_A\|\hat{R}^{-1}e$ with $c\asymp 1$ and $\|e\|\lesssim 1$, and we can assume $\|\hat{a}\|\asymp 1$ since $y_i$ are good approximation of $(A^\top A)^{-1}r_A$, then
\begin{align*}
    \hat{x}&=Y\hat{a} + h_2\\
    &=(A^\top A)^{-1}(r_A+h_1-(A^\top E_1+E_2)\hat{a})+h_2\\
    &=(A^\top A)^{-1}r_A + (A^\top A)^{-1}(h_1-E_2\hat{a})+\hat{R}^{-1}(\hat{R}h_2+\hat{R}(A^\top A)^{-1}A^\top E_1\hat{a}).\\
\end{align*}}

 {With $\|r_A\| = \|A^\top Ax^\ast\|\leq \|Ax^\ast\|\asymp \|\hat{R}^{-\top}r_A\|$, $\|Y\|\leq \sqrt{k}\max_i \|y_i\|\lesssim \sqrt{k}(\|x^\ast\|+c\kappa \|\hat{R}^{-\top}r_A\|)\lesssim \sqrt{k}\kappa \|Ax^\ast\|$ and $\|AY\|\leq \sqrt{k}(\|Ax^\ast\|+c\|\hat{R}^{-\top}r_A\|)$ we have following bounds
\begin{align*}
    \|h_1\|&\lesssim uk^2\|Ax^\ast\|,\\
    \|E_2\hat{a}\|&\leq un^\frac{3}{2}k\|Ax^\ast\|,\\
    \|h_2\|&\lesssim u\kappa k^2\|Ax^\ast\|,\\
    \|\hat{R}(A^\top A)^{-1}A^\top E_1 \hat{a}\|&\lesssim \|E_1\|\lesssim un^\frac{3}{2}\kappa \sqrt{k} \|Ax^\ast\|.
\end{align*}
Since $k$ is small, the result has the form
\begin{align*}
    x &= (A^\top A)^{-1}r_A+u\kappa n^{\frac{3}{2}}\|Ax^\ast\|\hat{R}^{-1}e_1+ un^\frac{3}{2}\|Ax^\ast\| (A^\top A)^{-1}e_2, \quad \|e_{1,2}\|\lesssim 1.
\end{align*}
The result consequently satisfies the condition of theorem \ref{iterativeIDS} as $u\kappa n^\frac{3}{2}\lesssim c$, thus the k-step Krylov-based SIRR solver is backward stable.
}

\end{document}